\documentclass[10pt]{article}
\usepackage[spanish,english]{babel}
\usepackage[utf8]{inputenc}
\usepackage{bbm}
\usepackage{dsfont}
\usepackage{amsmath}
\usepackage{amsthm}
\usepackage{amssymb}
\usepackage{amsfonts}
\usepackage{euscript}
\usepackage{amscd}
\usepackage{graphicx}
\usepackage{graphics}
\usepackage{epsfig}
\usepackage{latexsym}
\usepackage{fancyhdr}
\usepackage[all,cmtip]{xy}
\usepackage{srcltx}
\usepackage{titletoc}
\usepackage[bf]{titlesec}
\usepackage{mathrsfs}
\usepackage{mathabx}
\usepackage{verbatim}
\usepackage[T1]{fontenc}
\usepackage{ae,aecompl}
\usepackage{textcomp}
\usepackage{enumerate}
\usepackage{enumitem}
\usepackage[pdftex,linktocpage,breaklinks=true,colorlinks=true]{hyperref}
\usepackage{makeidx}
\usepackage{multirow}
\usepackage{multicol}
\usepackage{ifsym}
\usepackage{appendix}
\usepackage{calligra} 
\usepackage{textcomp}
\usepackage{color}
\usepackage[T1]{fontenc}
\usepackage{enumitem}
\setlist[enumerate]{
wide, 
nosep, 
labelwidth=*,
labelindent=0cm,
leftmargin=0cm}

\newtheorem{theorem}{Theorem}
\newtheorem{definition}{Definition}
\newtheorem{lemma}{Lemma}
\newtheorem{corollary}{Corollary}
\newtheorem{proposition}{Proposition}
\newtheorem{example}{Example}
\newtheorem{remark}{Remark}

\newcommand{\K}{\mathbb{K}}
\newcommand{\R}{\mathbb{R}}
\newcommand{\C}{\mathbb{C}}

\newcommand{\aff}{\hbox{\it Aff}\hskip.05cm}

\newcommand{\cp}{\mathbb{CP}}
\newcommand{\del}[2]{\frac{\partial #1}{\partial #2}}
\newcommand{\ZZ}{\mathbb{Z}}

\providecommand{\msc}[1]{\textbf{\textit{MSC:}} #1}
\providecommand{\keywords}[1]{\textbf{\textit{Keywords:}} #1}

\usepackage{graphicx}
\makeatletter
\newcommand*\bigcdot{\mathpalette\bigcdot@{.5}}
\newcommand*\bigcdot@[2]{\mathbin{\vcenter{\hbox{\scalebox{#2}{$\m@th#1\bullet$}}}}}
\makeatother

\begin{document}
	
\title{Plane polynomials and Hamiltonian vector fields 
\\
determined by their singular points}

\author{
 John A. Arredondo     $^a$ 
and
Jes\'us Muci\~no--Raymundo $^b$ 
\\
\small{$^a$ Fundaci\'on Universitaria Konrad Lorenz,
Bogot\'a, CP. 110231, Colombia,}
\\
\small{alexander.arredondo@konradlorenz.edu.co}
\\
\small{$^b$ Centro de Ciencias Matem\'aticas, 
UNAM, Campus Morelia, 
Mich. CP. 58089, M\'exico,}
\\
\small{muciray@matmor.unam.mx}
\\
} 

\date{\today} 

\maketitle


\begin{abstract} 
Let 
$\Sigma(f)$ be critical points of 
a polynomial $f \in \K [x,y]$
in the plane $\K^2$, 
where $\K$ is $\R$ or $\C$. 
Our goal is to study the 
critical point map $\mathfrak{S}_d$, 
by sending polynomials $f$ of degree $d$ to their
critical points $\Sigma(f)$.
Very roughly speaking,
a polynomial $f$ is essentially determined 
when any other  $g$ sharing 
the critical points of $f$ 
satisfies that $f = \lambda g$; 
here both are polynomials of at most degree $d$,
$\lambda \in \K^* $.
In order to describe the degree $d$ 
essentially determined polynomials,
a computation of the required number of 
isolated critical points 
$\delta(d)$ is provided. 
A dichotomy appears for the values of $\delta(d)$; 
depending on a certain parity
the space of essentially determined polynomials
is an open or closed Zariski set.  
We compute the map $\mathfrak{S}_{3}$, 
describing under what conditions 
a configuration of four points 
leads 
to a
degree three essentially determined polynomial.
Furthermore, we describe explicitly 
configurations supporting
degree three 
non essential determined
polynomials. 
The quotient space of essentially 
determined polynomials of degree three up to 
the action of the affine group $\aff(\K^2)$
determines a singular surface over $\K$. 
\end{abstract}

\msc{35B45; 35R09; 35B65; 35B33}


\section{Introduction}

Let $\K = \R$ or $\C$. We then ask
under what conditions a 
polynomial $f \in \mathbb{K}[x,y]$ is 
essentially determined
by its critical points $\Sigma (f) \subset \K^2$?
Thus, we want to study the 
{\it critical point map} sending 
polynomials of degree $d$ 
to their critical points
\begin{equation}\label{mapeo-de-puntos-criticos}
\mathfrak{S}_d : f \longmapsto \Sigma(f) 
\, ,
\end{equation}
\noindent 
where $\Sigma(f) \doteq I(f_x, f_y)$ 
is the affine algebraic variety
(not necessarily reduced)
generated by the ideal 
of partial derivatives of $f$, see
Definition \ref{mapeo-de-puntos-criticos-I}.
Our approximation route uses 
a finite dimensional framework. Let
$\K[x,y]^0_{\leq d}$
be the $\K$--vector space of polynomials 
having 
at most degree $d$ ($\geq 3$) 
and zero independent term, 
and let $\mathcal{P}=\{ (x_\iota, y_\iota) \}$ 
be a configuration of $n$ different  points in the plane.      
The linear projective subspace of the polynomials
with 
critical points at least in $\mathcal{P}$, denoted as
\begin{equation} 
\label{un-espacio-de-polinomios}
\mathcal{L}_{d}(\mathcal{P}) \doteq 
Proj \big(
\{ f \in \K[x,y]^0_{\leq d } \ \vert \
\mathcal{P}  \subseteq \Sigma (f)\} \big),
\end{equation}
 
\noindent 
is well defined.  
We say 
that a polynomial
$f$ is {\it essentially determined }
by $\mathcal{P} $ 
when $\mathcal{L}_d(\mathcal{P})$ is 
a projective point
$\{ \lambda f \ \vert \ \lambda \in \K^*\}$,
see Definition \ref{definicion-clave}.
All this leads us to the following.

\medskip

\noindent
\textbf{Interpolation problem for critical points.}
Let 
$\mathcal{P} \subset \K^2$ be a configuration of $n$ 
different points, 
we try to
determine the projec\-tive subspace 
$\mathcal{L}_{d}(\mathcal{P}) $
of polynomials of at most degree $d$ 
with critical points at least 
in
$\mathcal{P}$.

\medskip

This problem has several novel features. 
The critical values 
$\{ c_\iota \} \subset \K$ of $f$
can appear in different level curves 
$\{ f(x,y)-c_\iota=0  \}$; 
it is natural in Hamiltonian vector field theory 
and moduli spaces of polynomials, see 
P. G. Wightwick \cite{Wightwick}
and
J. Fern\'andez de Bobadilla \cite{Fernandez-de-Bobadilla}. 
This is a main difference with the widely considered 
problem of linear system of curves 
in $\cp ^2$, {\it e.g.} 
R. Miranda, 
\cite{Miranda}
and
C. Ciliberto \cite{Ciliberto}. 

Very roughly speaking, for degree $d \geq 3$  
the relevant data are the
cardinality and position of the configuration 
$\mathcal{P}$, as candidate to be a
critical point configuration $\Sigma(f)$.
For degree three, the 
prescription of four critical points 
is suitable.  
For degree $d \geq 4$, 
however, 
the generic configuration $\mathcal{P}$
having $(d-1)^2$ points is too restrictive. 
Thus the fiber $\mathfrak{S}_d^{-1}(\mathcal{P})$ 
will be generically empty. 
It follows that, the position of the configurations 
$\mathcal{P}$ coming from polynomials 
is the hardest part to be characterized. 
At this first stage, we consider mainly 
$\mathcal{P}$ as isolated points of multiplicity one,
Remark 
\ref{remark-de-multiplicidades} provides an 
explanation. 
Our first result describes the role of cardinality 
$\delta(d)$ of $\mathcal{P}$
in Eq. \eqref{un-espacio-de-polinomios}, 
see Proposition \ref{abierto-cerrado-Zariski}.
\smallskip

{\it 
Dichotomy on the required number of critical points

\noindent 
If 
the dimension of  $\K[x,y]^0_{\leq d}$ is odd 
(resp. even) 
then the configurations $\{\mathcal{P} \}$ with
$\delta(d)$  points
and $dim_\K(\mathcal{L}_d(\mathcal{P})) \geq 0$
determine an open
(resp. closed) 
Zariski set in the space
of configurations with $\delta(d)$ points, denoted as
${\it Conf}(\K^2, \delta(d))$. 
}

\smallskip 

We compute the 
critical point map $\mathfrak{S}_{3}$. Thus,
a description for the four critical point configurations
$\{ \mathcal{P}\}$ with essentially determined polynomials is
provided.
Recall that the affine group $\aff(\K^2)$
acts on the space of polynomials, 
see Eq. \eqref{accion-afin}.
This action is rich enough and yet treatable for 
degree three. 
Let 

\quad

\centerline{$\mathscr{A} \doteq 
\{x_4 y_4(x_4 + y_4 -1)(x_4 + y_4)(x_4-1)
(y_4 -1)=0\}
\subset \K^2= \{(x_4, y_4)\}
$}

\quad

\noindent 
be an arrangement of six lines 
from two nested  triangles, 
one of them is 
$\triangle=\{ (0,0), (1, 0), (0,1) \}$;
see Fig. \ref{mos}.a.
We prove the following result.

\begin{theorem}\label{segundo-teorema-grado-tres}
Let $f$ 
be a degree three polynomial  
having  at least four 
critical points $\Sigma(f)$. 

\smallskip

\begin{enumerate}
\item  
$f$ is essentially determined if and only if 
up to affine transformation the four points are
$$
\Sigma(f) = \{ (0,0), (1,0), (0,1), (x_4, y_4) \}
\hbox{\  and \ }
(x_4, y_4) \notin \mathscr{A}.
$$

\item
$f$ is not essentially determined if and only if
up to affine transformation the four points are
$$
\{ (0,0), (1,0), (0,1), (x_4, y_4) \} 
\hbox{\  and \ }
(x_4, y_4) \in \mathscr{A}. 
$$ 
Moreover, in this case
$\Sigma(f)$ can be four isolated points 
or two parallel lines.
\end{enumerate}
\end{theorem}

\smallskip 
In simple words,  
the 4--th point $(x_4, y_4)$ 
generically determines the polynomial $f$. 
We compute the {\it fundamental domain} for this $\aff(\K^2)$--action,  obtaining 
a tessellation of $\K^2= \{(x_4, y_4)\}$ with 24 tiles, 
a seen in 
Fig. \ref{transformaciones-gauge-2}.

As is expected, some interesting phenomena 
occur for configurations with non trivial 
isotropy groups. 
For degree $d \geq 3$, a particular family of 
configurations is the  
grids of $(d-1)^2$ points from the intersection of 
two families of $d$ parallel lines in $\K^2$,
see
Definition \ref{red-d2-puntos}. 
They provide examples of 
non essential determined 
polynomials with $(d-1)^2$ 
Morse critical points.
A remaining open question, are these 
grids of $(d-1)^2$ points
the unique mechanism in order to produce
non essential determined Morse polynomials?

From the point of view of vector fields, 
we are studying 
under what conditions the zeros a 
Hamiltonian vector field
determine it in a unique way? 
This is a very general and 
interesting issue in real and complex
foliation theory, studied by; 
X. G\'omez--Mont, G. Kempf \cite{GomezMont-Kempf}, 
J. Artes, J. Llibre, N. Vulpe \cite{Artes-Llibre-Vulpe},
A. Campillo, J. Olivares \cite{Campillo-Olivares} 
and V. Ram\'irez \cite{Ramirez} see Corollary 
\ref{espectros}.
Related works are accurately described 
in Section \S \ref{closing-remarks}.

The content of this work is as follows.
In \S2--3,  we study the problem of the dimension of 
linear systems for polynomials with critical points, 
using the degree as parameter.
In section \S4, 
we characterize  polynomials essentially determined by 
their configurations  
of critical points; this proves Theorem 1.
In section \S5, 
we focus  in the degree four case.
For each configuration of six points, we obtain
a plane curve of degree six parametrizing
the essentially 
determined polynomials, see 
Proposition \ref{curvas-de-interpolacion}.
Section \S6 explores the behavior of
pencils of Hamiltonian vector fields with 
common simple zeros.


\section{Linear systems $\mathcal{L}_{d} (\mathcal{P})$   }
Let
$\K[x,y]^0_{\leq d}$ (resp. $\K[x,y]^0_{= d}$)
be the $\K$--vector space of polynomials 
having at most degree $d \geq 3$ (resp. the set for degree $=d$) 
and zero independent term.
Consider
\begin{equation}\label{polinomio-general-grado-d}
f(x,y)=\sum_{
1 \leq \iota +j \leq d} 
a_{\iota j} x^\iota y^j 
\ \
\in \mathbb{K}[x,y]_{\leq d}^0,
\end{equation}

\noindent 
from which
the $\K$--dimension of $\K[x,y]_{\leq d}^0$ is
$\frac{1}{2}(d^2 + 3d)$, and
its projectivization is
\begin{equation}
\label{proyectivizacion}
Proj\big( \K[x,y]^0_{\leq d} \big) 
= 
\big\{ [f] \ \vert \ f \in \K[x,y]^0_{\leq d} \big\} 
=
\K \mathbb{P}^{\frac{1}{2}(d^2 + 3d - 2)},
\end{equation}

\noindent
where $[ \ \ ]$ denotes a projective class. 
Recall that 
\begin{equation}
{\it Conf}(\K^2,n) = 
\big\{ \, 
\mathcal{P}=
\{ (x_1, y_1), \ldots , (x_n, y_n )\} 
\ \vert \ 
(x_\iota, y_\iota) \neq (x_j , y_j ) \hbox{ for } \iota \neq j
\big\} / Sym (n)
\end{equation}

\noindent 
is the {\it space of unordered configurations of $n$ points}
in $\K^2$, where
the symmetric group in $n$ elements, $Sym(n)$,
acts by exchanging the points.  The configuration space 
${\it Conf}(\K^2, n)$ is a $\K$--analytic manifold.

\begin{definition} 
\begin{upshape}
Given a configuration 
$
\mathcal{P}
\in 
{\it Conf}(\K^2,n ),
$
the {\it 
linear system of 
polynomials of at most degree $d$ 
with critical points at least 
in $\mathcal{P}$}
is the projective subspace
\begin{equation}
\label{contencion-basica}
\mathcal{L}_{d}(\mathcal{P}) = 
\big\{
[f] \ \vert \ 
\mathcal{P} \subseteq 
\{f_x(x,y) =0\} \cap \{ f_y(x,y) =0 \} 
\big\} 
\subset 
Proj\big( \K[x,y]^0_{\leq d} \big).
\end{equation}
\end{upshape}
\end{definition}

In algebraic geometry language,
$\{ f_x (x,y) = 0\}$, $\{f_y (x,y) =0\}$ 
belong to the linear system of 
algebraic curves 

\centerline{
$\mathcal{L}_{d-1}
\big(
- \Sigma_{\alpha = 1}^n (x_\iota , y_\iota ) 
\big)$ }

\noindent
see \cite{Miranda} and \cite{Ciliberto}. 
In several places however, 
we consider $f_x$, $f_y$ as functions 
and not just as algebraic curves.

The polynomials of at most degree $d$, the  
polynomial Hamiltonian vector fields 
and the polynomial vector fields, of at most degree $d-1$, 
are related by linear maps
$$
\begin{array}{ccccc}
\K[x,y]^0_{\leq d}  &  \stackrel{\cong}{\longleftrightarrow} 
&  Ham(\K^2)_{\leq d-1} 
& \longrightarrow & \mathfrak{X} (\K^2)_{\leq d-1} 
\\
&&&& \vspace{-.3cm}
\\
f & \longleftrightarrow & X_f = -f_y \del{}{x} + f_x \del{}{y} & \longrightarrow & X_f.
\end{array}
$$

\noindent 
In the space of Hamiltonian vector fields,
$\mathcal{L}_{d} (\mathcal{P})$ 
determines  a linear subspace
$$
\{  \lambda X_f \ \vert \
\mathcal{P}  \subseteq 
\mathcal{Z}(\lambda X_f),\ \lambda \in \K^*
\} \subset  Ham(\K^2)_{\leq d-1}; 
$$
set theoretically, the zeros 
$\mathcal{Z}(\lambda X_f)$ of the vector field $X_f$
coincide
with $\{f_x(x,y)=0 \} \cap \{f_y(x,y)= 0\}$.

\begin{definition}
\label{definiciones-de-multiplicidad}
\begin{upshape}
Let $f \in \K[x,y]$ be a non constant polynomial. 
Over $\K=\C$, the {\it Milnor number of $X_f$ at a 
zero point $(x_\iota, y_\iota) \in \mathcal{Z}(X)$}
is 

\centerline{$ 
\mu_{(x_\iota, y_\iota)}(X) =
\dim_\C \frac{\mathcal{O}_{\C^2, (x_\iota, y_\iota)}  
}{ <-f_{y}, f_x  >}
$,}

\noindent 
where $\mathcal{O}_{\C^2, (x_\iota, y_\iota)}$ is the local ring
of holomorphic functions at the point $(x_\iota, y_\iota)$
and $<-f_{y}, f_x  >$ is the ring generated by the partial derivatives. 
\end{upshape}
\end{definition}

\begin{remark}
\label{remark-de-multiplicidades}
\begin{upshape}
1. 
Over $\K=\C$, if $(x_\iota, y_\iota)$ is
an isolated singular point of $f$, 
then the notions  
of multiplicity for the intersection of the curves
$\{f_x(x,y)=0\} \cup \{f_y(x,y)=0\}$  
and the Milnor number for $X_f$ coincide;
see \cite{Jong-Pfister} p.\,174.

\noindent 2. A priori,  
we consider each point $(x_\iota, y_\iota) \in \mathcal{P}$
in \eqref{contencion-basica}
with multiplicity of intersection one for the algebraic curves 
$\{ f_x(x,y)= 0\}$ and $\{ f_y(x,y)=0\}$.

\noindent 
3. By B\'ezout's theorem, 
the maximal number of isolated singularities 
of $X_f$ on $\C^2$
is $(d-1)^2$. 
In this case   
all the affine singularities are of multiplicity one. 

\noindent 
4. Moreover, the maximal number of isolated singularities 
of $X_f$ extended to $\mathbb{CP}^2$ is

\centerline{$
(d-1)^2 + d
$.}

\noindent 
Here the upper bound $d$ comes from the intersection of a 
generic projectivized level curve 
$\{ f=c\}$ with the line at infinity;
see
\cite{GomezMont-Kempf}, \cite{Campillo-Olivares} for the
case of rational vector fields, 
which are not necessarily Hamiltonian.
\end{upshape}
\end{remark}

Let $\mathbb{A}^2_\K = \hbox{Spec}\, \K[x,y]$ be the affine 
scheme of the affine plane $\K^2$, see 
\cite{Eisenbud-Harris} pp. 48--49.

\begin{definition}
\begin{upshape}
The {\it critical point map of degree $d$} is the map
\begin{equation}
\label{mapeo-de-puntos-criticos-I}
\begin{array}{rcl}
\mathfrak{S}_d 
: \K[x,y]_{= d}
&\longrightarrow &  \hbox{Spec}\, \K[x,y]
\\
&& \vspace{-.2cm}
\\
f &\longmapsto &\Sigma(f)= I(f_x, f_y) \, ,
\end{array}
\end{equation}
\noindent 
sending a
polynomial of degree $d$ 
to its critical points $\Sigma(f)$, 
as an affine algebraic variety
(not necessarily reduced) generated by the ideal 
of partial derivatives of $f$. 
\end{upshape}
\end{definition}

In fact, $\Sigma(f)$ can be understood as a subscheme, 
with support at the points  
$\{f_x(x,y)=0\}\cap \{f_y(x,y)= 0\}$, where the 
sheaf of ideals is defined by the germs of $I(f_x, f_y)$;
compare with 
\cite{Campillo-Olivares},
\cite{Eisenbud-Harris} p.\,100.
In a set theoretical language,   
$\Sigma(f)$ determines points and even algebraic
curves.
However  
in the study of rational vector fields on $\mathbb{CP}^2$, the case
of foliations having singularities along curves
is removed, see
\cite{GomezMont-Kempf}, \cite{Campillo-Olivares}. 

\begin{remark}
\begin{upshape}
The simplest case of the 
interpolation problem for singular points
occurs when $\Sigma(f)$ is a finite 
set of points of multiplicity one, {\it i.e.}
$\{ f_x(x,y)=0\}$ and $\{f_y(x,y) = 0\}$ 
have transversal intersections. The
$\Sigma(f)$ is a configuration in 
${\it Conf}(\K^2, n)$, for $0 \leq n \leq (d-1)^2$.
\end{upshape}
\end{remark}

Our former task is as follows:
{\it Given a 
configuration $\mathcal{P}$,
which is $dim_\K (\mathcal{L}_{d}(\mathcal{P}))$?}

\noindent 
To be clear,
three relevant data must be considered 
the degree $d$ of the polynomials $\{ f \}$, 
the cardinality $n$ and the position
of the configuration $\mathcal{P}$. The following 
diagram explains: 

\begin{center}
\begin{picture}(400,75)
\put(-40,5)
{\vbox{\begin{equation}\label{diagrama-maestro}
\end{equation}}}
\put(-7,4){position of $\mathcal{P}$}
\put(-28,46){cardinality $n$ of $\mathcal{P}$}

\put(60,7){\vector(2,1){30}}
\put(60,49){\vector(2,-1){30}}

\put(90,25){
$\begin{array}{c}
dim_\K (\mathcal{L}_{d}(\mathcal{P}))
\end{array}=
$}

\put(186,25){
$\left\{ \begin{array}{lcl}
 -1  &  & \mathcal{L}_{d}(\mathcal{P}) = \emptyset. \\
& & \vspace{-.1cm} \\
 \ 0 & & [f] = \mathcal{L}_{d} (\mathcal{P}) 
= \K \mathbb{P}^0   \\
&& f \hbox{ is essentially determined.}\\
& & \vspace{-.1cm} \\
\kappa \geq 1 & & 
[f] \in \mathcal{L}_{d} (\mathcal{P}) 
= \K \mathbb{P}^\kappa 
\\
& & f \hbox{ is non essential determined.}
\end{array}\right.
$}
\end{picture}

\vskip1cm
\end{center}

\noindent
The natural concepts are as follows.

\begin{definition}\label{definicion-clave} 
\begin{upshape}
Let $f \in \K[x,y]_{\leq d}^0$ be a polynomial and
let $\mathcal{P}$ be a
configuration of $n$ points in $\K^2$.
\begin{enumerate}[leftmargin=*]
\item
A polynomial $f$ is {\it essentially determined by}
$\mathcal{P}$  when 
$ [f]=\mathcal{L}_{d}(\mathcal{P})  $.

\item
A polynomial $f$ is {\it non essential determined
by}
$\mathcal{P}$  when 
$ [f]\in \mathcal{L}_{d}(\mathcal{P})  $ and
$dim_\K( \mathcal{L}_{d}(\mathcal{P})) \geq 1$.

\item 
$\mathcal{P}$ is {\it a forbidden configuration}
(for polynomials of at most degree $d$) 
when 
$\mathcal{L}_{d}(\mathcal{P})= \emptyset$. 

\item
The {\it set of degree $d$ essentially
determined polynomials} is
\begin{equation}
\label{todos-los-polinomios-esencialmente-det}
\mathcal{E}_{d}  \doteq
\bigcup _{ \tiny{\mathcal{P} }  }
\mathcal{L}_{d} (\mathcal{P})
\subset 
Proj \big(\K[x,y]^0_{\leq d} \big),  
\end{equation}
where the union is over all 
configurations $\{ \mathcal{P} \}$ such that 
$dim_\K(\mathcal{L}_{d} (\mathcal{P}))= 0$.
\end{enumerate} 
\end{upshape}
\end{definition}

\begin{remark}
\begin{upshape}
\begin{enumerate}
\item 
The strict set theoretical inclusion $\mathcal{P} \varsubsetneq \Sigma(f)  $
can be satisfied for essentially determined polynomials $f$, 
for example as with the case of a product of three lines
one with multiplicity two, say $f=L_1^2 L_2$.

\item 
The set of degree three essentially determined polynomials 
$\mathcal{E}_{3}$ is a union of projective spaces, however 
it is not a projective space, as Proposition 
\ref{segundo-teorema-grado-tres} will show. 

\item
As is expected,
many of the projective classes 
in $\mathcal{E}_{d}$ arise from Morse polynomials. 
The converse is not true, see  Corollary 
\ref{polinomios-Morse-no-esencialmente-determinados}. 
\end{enumerate} 
\end{upshape}
\end{remark}


\section{On the number of required critical points}
\label{seccion-polinomios-genericos}

A novel aspect of the interpolation problem for
critical points is  its cardinality; 
the configurations having a certain number 
$\delta(d)$ of points 
determine open or closed Zariski sets in 
$\K[x,y]_{\leq d}^0$. 
As a key point, the dimension 
$\frac{1}{2} (d^3+ 3d)$ of $\K[x,y]_{\leq d}^0$
can be even or odd. Starting with degree $d=4$, 
the pattern of these dimensions is 4--periodic; 
even, even, odd odd, $\ldots$ .
See the third column in Table \ref{delta}.

\begin{proposition}\label{abierto-cerrado-Zariski}
(A dichotomy on the number $\delta(d)$ of required critical points)
\
Let $ \K[x,y]^0_{\leq d}$
be the set of polynomials having at most degree $d \geq 3$
and let
\begin{equation}
\label{valores-de-delta-explicitos}
\delta (d) 
\doteq  
 \left\{ \begin{array}{cl}
\frac{1}{4}(d^2+3d-2)& \hbox{ when } \frac{1}{2}(d^2+3d) \hbox{ is odd}, \\
& \vspace{-.2cm}\\
\frac{1}{4}(d^2+3d) & \hbox{ when } \frac{1}{2}(d^2+3d) \hbox{ is even} .             
\end{array}\right. 
\end{equation}

\noindent 
1. If 
the dimension of  $\K[x,y]^0_{\leq d}$ is odd, 
then the configurations $\{\mathcal{P} \}$ with
$\delta(d)$  points
and $dim_\K(\mathcal{L}_d(\mathcal{P})) \geq 0$
determine an open Zariski set in ${\it Conf}(\K^2, \delta(d))$. 

\smallskip 

\noindent 
2. If 
the dimension of  $\K[x,y]^0_{\leq d}$ is even, 
then the configurations $\{\mathcal{P} \}$ with
$\delta(d)$ points
and $dim_\K(\mathcal{L}_d(\mathcal{P})) \geq 0$
determine a closed Zariski set in ${\it Conf}(\K^2, \delta(d))$.
\end{proposition}

\begin{proof}
Let $ f(x,y) \in \K[x,y]^0_{\leq d}$ be a polynomial 
as in \eqref{polinomio-general-grado-d}.
Assume that $\mathcal{P}=\{ (x_\iota, y_\iota) 
\ | \  \iota =1,\ldots, n\} $ 
is set theoretically contained in $\Sigma(f)$.
A priori,
each point $(x_\iota , y_\iota ) \in \mathcal{P}$ 
will drop the dimension of the vector space 
$\K[x,y]^0_{\leq d} $ 
by two. 
In the linear framework,
this leads to a linear system of 
$2n$ equations:
\begin{equation}
\label{critical-point-equations}
f_x(x_\iota,y_\iota) = f_y(x_\iota,y_\iota) =0,
\ \ \ \ \iota =1,\ldots, n,
\end{equation}
with  $\{a_{\iota j}\}$ as variables.
Following B\'ezout's theorem for a moment, 
let us consider a configuration with 
$n=(d-1)^2$ points. 
We have a linear map
\begin{equation}\label{matriz-phi}
\begin{array}{rcl}
\phi: \K [x,y]^0_{ \leq d} \cong 
\K^{\frac{1}{2}(d^2 + 3d)} & \longrightarrow & 
\K^{2(d-1)^2} \\
f & \longmapsto &
\big(f_x(x_1, y_1), \ldots , f_x(x_{(d-1)^2}, y_{(d-1)^2}), \,
f_y(x_1, y_1), \ldots , f_y(x_{(d-1)^2}, y_{(d-1)^2})\big).
\end{array}
\end{equation}
\noindent 
The interpolation matrix 
$\phi$ depends on $\mathcal{P}$, 
and by notational simplicity we omit this dependence.  
The matrix $\phi$ has 
$\frac{1}{2}(d^2 + 3d)$ columns,
$2 (d-1)^2$  rows and a very particular shape 
because of the partial derivatives involved  in it,
see Eqs. 
(\ref{matriz-grado-tres}), 
(\ref{matriz-grado-cuatro})
for explicit examples with $d= 3$, $4$.

For degree $d = 3$ 
and a configuration $\mathcal{P}$ of 4 points;
however then the rank of the matrix 
$\phi$ associated to $\mathcal{P}$ is 8
if and only if    
$dim_\K (\mathcal{L}_{3} (\mathcal{P}) )= 0$. 
If we consider  
degree $d \geq 4$, then  
the number of rows of $\phi$ 
is bigger than the number of columns.
We must reduce the number $n$ of required points in the configurations 
$\mathcal{P}$, this $ n < (d-1)^2$.
The number $\delta(d)$ in 
\eqref{valores-de-delta-explicitos} determines two 
possibilities.

\smallskip

\noindent 
{\it Case 1 in \eqref{valores-de-delta-explicitos}.} 
For $\mathcal{P}$ with
$\delta(d)=\frac{1}{4}(d^2+3d-2)$ points,
the interpolation matrix  
$\phi$ has   $\frac{1}{2}(d^2+3d)$ odd columns and 
$\frac{1}{2}(d^2+3d-2)$ even rows, for example for 
$(d+1)=3, 6, 7$. 
Moreover, 

\centerline{
(number of columns of $\phi$)--1 = 
(number  of rows of $\phi$).}

\noindent 
The dimension of the kernel of $\phi$ is at least one, 
thus $dim_K(\mathcal{L}_d(\mathcal{P})) \geq 0$. 
There are $\frac{1}{2}(d^2+3d)$ minors $A_j$ from the matrix
$\phi(x_1, y_1, \ldots , x_{\delta(d)}, y_{\delta(d)})$.
The complement of the algebraic equations

\centerline{
$\{ \Pi_j det(A_j(x_1, y_1, \ldots , x_{\delta(d)}, y_{\delta(d)}))  =0 \} \subset {\it Conf}(\K, \delta(d))$}

\noindent 
describes the set of configurations having
$dim_K (\mathcal{L}_d(\mathcal{P}))=0$, 
corresponding to 
the essentially determined polynomials. 
These configurations of $\delta(d)$ points in 
${\it Conf}_{\delta (d)} (\K^2)$
determine an 
open Zariski and dense set, that is the second part of 
assertion (1).

\smallskip

\noindent 
{\it Case 2 in \eqref{valores-de-delta-explicitos}.} 
The dimension of $\K[x,y]_{\leq n}^0$ is even
and we assume  $\frac{1}{4}(d^2+3d) \in \mathbb{N}$
points in $\mathcal{P}$.
The interpolation matrix  
$\phi$ is square of even size, and 
there are $\frac{1}{2}(d^2+3d)$ 
columns and rows; for example when $d= 4, 5$.

\noindent 
If we assume $\mathcal{P}$ such 
that $\{ det(\phi(x_1, y_1, \ldots , x_{\delta(d)},y_{\delta(d)}))\neq 0 \}$, 
then the only vector in the $\{ a_{\iota j} \}$ variables 
solving the linear system 
\eqref{critical-point-equations} is zero.
The set of desired polynomials is empty.

\noindent 
The configuration with non empty polynomials 

\centerline{$\left\{
\mathcal{P} \ \vert \ det(\phi(x_1, y_1, \ldots , x_{\delta(d)},y_{\delta(d)}))\neq 0 \right\}
\subset {\it Conf}(\K, \delta(d))
$} 

\noindent
determines an algebraic set.
\end{proof}

\begin{table}[htbp]
\begin{center}
\begin{tabular}{|c|c|c|c|c|}
\hline
&&&& \vspace{-.4cm}
\\
& & number of & number of & 
Zariski
\\
degree $d$&$\delta(d)$&   
columns in $\phi$ & rows in $\phi$ & topology of
\\
&eq.  
\eqref{valores-de-delta-explicitos}& $\frac{1}{2}(d^2+3d)$ & $2 \delta(d)$ &
$\{\mathcal{P} \} \subset {\it Conf}(\K^2, \delta(d))$
\\
&&&& \vspace{-.3cm}\\
\hline
\hline
3 & 4&  9 & 8 & closed\\
\hline
4 & 7&  14 & 14 & open \\
\hline
5 & 10 &  20 & 20 & open \\
\hline
6  & 13&  27 &  26 & closed\\
\hline
7  & 17&  35 &  34 & closed \\
\hline
\end{tabular}
\caption{Dimensions and values 
for the interpolation problem.}
\label{delta}
\end{center}
\end{table}
\bigskip

Recalling \eqref{proyectivizacion}, 
the {\it expected  projective dimension of
$ \mathcal{L}_{d}(\mathcal{P})$}, which is
the
linear system of polynomials of at most degree $d$
with critical points at least in 
$\mathcal{P} \in {\it Conf}(\K^2, n)$, 
is

\centerline{$max
\left\{ 
\frac{1}{2}(d^2+3d-2)-2n , \ -1 \right\} 
$. 
} 

\noindent 
In Section \ref{polinomios-grado-cuatro}, 
we provide an alternative for studying 
the even dimension case in 
Proposition \ref{abierto-cerrado-Zariski}.

\section{Essentially determined polynomials of degree three}
\label{sec.-conf.-degree-3}

\subsection{A linear system}
 
\begin{figure}[h!]
\centerline
{\includegraphics[scale=0.58]{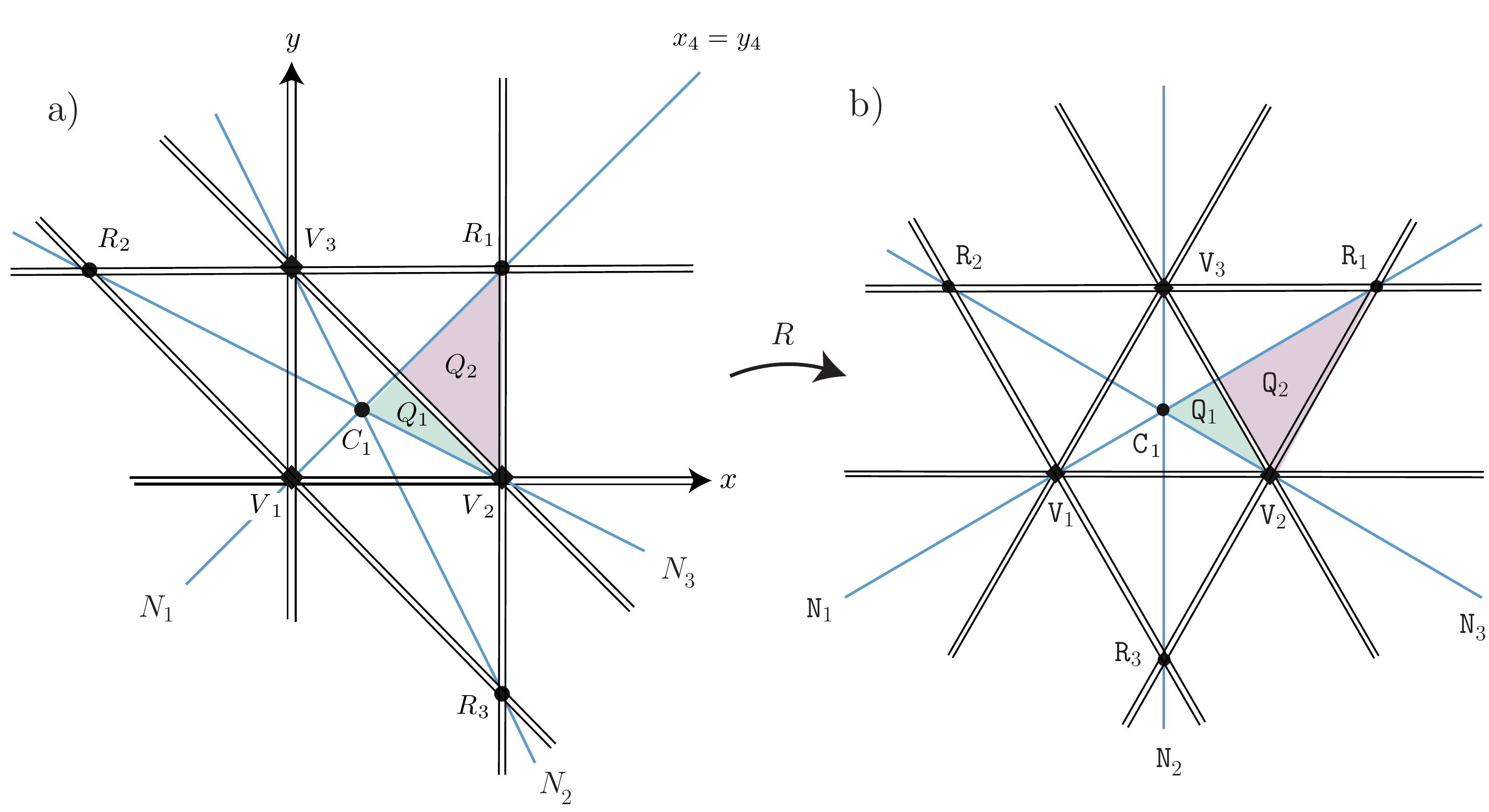}}
\caption{ 
a) The line arrangement $\mathscr{A}$
(of double lines) and the triangle
$\triangle= \{V_1, \, V_2, \, V_3\}$.
b)  
The analogous objects under 
the linear map $R$, sending   
$\mathscr{A}$ to ${\tt A}$ and
$\triangle$ to $\Delta$. 
\label{mos}
}
\end{figure}

In order to apply elementary methods,
we introduce a very simple configuration 
of four  points, depending essentially of 
the fourth one $(x_4, y_4)$. Secondly,
we must find a polynomial $f(x_4, y_4, x,y)$ with 
a critical point set containing  
the above simple configuration. Let 
\begin{equation}
\label{arreglo-de-seis-lineas}
\mathscr{A} \doteq
\big\{ x y (x+y-1) (x+y) (x-1) (y-1) = 0\big\}
\end{equation}
be an arrangement of six $\K$--lines; 
it is illustrated in Fig. \ref{mos}.a.

\begin{lemma}
Let 

\centerline{$\mathcal{P}= 
\{V_1=(0,0), \, V_2=(1,0), \, V_3 =(0,1), \, (x_4, y_4) \} 
\in {\it Conf}(\K^2, 4)$, 
\ \ \
$(x_4, y_4) \notin \mathscr{A}$,  
}

\noindent  
be a four point configuration. 
The polynomial 
\begin{align}
\label{solucion-general-grado-3-usando-base-canonica-0}
f(x_4, y_4, x, y)=&
(
y_4^2(y_4-1)(-1+2 x_4 +y_4)(2x^3 - 3x^2)+x_4^2(x_4-1) ( -1 + x_4+2 y_4)(2y^3  -  3y^2)
 \\
-&6x_4y_4(x_4-1)(y_4-1) (x^2y +xy^2-xy) 
)
a_6  \quad \in \K[x,y]_{=3},\nonumber
\end{align}
for $a_6 \in \K^*$,
is well defined and 
$\mathcal{P} =
\Sigma\big( f(x_4, y_4, x,y) \big)$. 
\end{lemma}

It shall be convenient to write the Eq.
\eqref{solucion-general-grado-3-usando-base-canonica-0}
as a map to the space of polynomials

\begin{equation}\label{solucion-general-grado-3-usando-base-canonica}
f(x_4, y_4, \  \ , \ \ ):
\K^2 \backslash \mathscr{A}  \longrightarrow  
\K[x,y]_{=3},
\ \ \
(x_4, y_4) \longmapsto 
f(x_4, y_4, x,y).
\end{equation}

\begin{proof} 
Let the following be a polynomial 
\begin{equation}\label{polinomio-grado-tres}
f(x,y)=a_1 x^3+a_2 x^2 y+a_3 x y^2+a_4 y^3+a_5 x^2
+a_6 x y+a_7 y^2+a_8 x+a_9 y 
\ \in \K[x,y]^0_{\leq 3}.
\end{equation}

\noindent
By notational simplicity, only one subindex
$a_{\iota}$ is considered.
Let 
$\{ (x_\iota, y_\iota) \ \vert \ 
\iota = 1,\ldots, 4\}$ be an arbitrary  
configuration, and we require 
$(a_1, \ldots, a_9 )$ to be solutions of the 
linear system

\begin{equation}\label{matriz-grado-tres}
\left(
\begin{smallmatrix}
&&& \vdots &&&&& \\
 3x_\iota^2 & 2x_\iota y_\iota & y_\iota^2 & 0 & 2x_\iota & y_\iota & 0 & 1 & 0 \\
 0 & x_\iota^2 & 2x_\iota y_\iota & 3y_\iota^2 & 0 & x_\iota & 2y_\iota & 0 & 1 \\
&&& \vdots &&&&&
\end{smallmatrix}
\right) 
\left(
\begin{smallmatrix}
a_1 \\
\vdots \\
a_9
\end{smallmatrix}
\right)
= 
\left(
\begin{smallmatrix}
0 \\
\vdots \\
0
\end{smallmatrix}
\right) .
\end{equation}

\noindent
The interpolation matrix
$\phi$ in \eqref{matriz-grado-tres} has
9 columns and 8 rows. The choice 
$\mathcal{P} = \{ (0,0), (1,0), (0,1), (x_4, y_4)\}$
determines the linear system with only two equations
$$
\begin{array}{c}
f_x (x,y) = 3a_1 x^2 + 2a_2 xy + a_3 y^2 + 2a_5 x + a_6 y + a_8 =0, \,
\\
f_y (x,y) = a_2 x^2 + 2a_3 xy + 3a_4 y^2 + a_6 x + 2a_7 y + a_9=0.
\end{array}
$$

\noindent
Obviously, $(0,0) \in \mathcal{P}$
implies the vanishing of the linear part 
$
f_x (0, 0) = a_8 = 0 = a_9 =f_y(0,0)$.
The linear conditions imposed by
$(1, 0)$ and $(0,1)$ are
$$
\left\{
\begin{array}{lccl}
f_x (1, 0) = 3a_1 + 2a_5 =0 &&& a_1 = - \frac{2}{3} a_5, \\
f_y (1, 0) = a_2 + a_6 =0 &&& a_6 = - a_2,
\end{array} \right.
$$
$$
\left\{ 
\begin{array}{lccl}
f_x (0, 1) = a_3 + a_6  = 0 & & & a_6= -a_3, \\
f_y (0, 1) = 3a_4 + 2a_7 = 0 & & & a_4 = \frac{-2}{3} a_7.
\end{array}  \right.
$$

\noindent
The solution of this system
\begin{equation}
\label{solucion-general-grado-3-usando-base-canonica-0-2}
\begin{array}{ll}
f(x_4, y_4, x, y)=&
a_6 \left(
\dfrac{ y_4 (-1+2 x_4 +y_4)}{3x_4 (x_4-1) }x^3 
- x^2y - xy^2 + 
\dfrac{x_4( -1 + x_4+2 y_4)}{3 y_4(y_4-1) }y^3  
\right.
\\
& \vspace{-.3cm}
\\
& \left. \ \  \ \
+ \dfrac{y_4(1-2 x_4 -y_4)}{2 x_4 (x_4-1) }x^2
+ xy + \dfrac{x_4(1-x_4 -2 y_4)}{2 y_4 (y_4-1) }y^2\right)
\  \
\in \K[x,y]_{=3}
\end{array}
\end{equation} 

\noindent  
has rational coefficients.
If we normalize, we get Eq. \eqref{solucion-general-grado-3-usando-base-canonica-0}.
\end{proof}

\begin{corollary}
Let 

\centerline{$\mathcal{P}_1=
\{(0,0), (1,0), (0,1), R_1 \doteq (1, 1) \} \in {\it Conf}(\K^2, 4)$}

\noindent 
be  a four point configuration,  
then $dim_\K (Proj(\mathcal{L}_3(\mathcal{P}_1 )))=1$.
\end{corollary}

We say that, $R_1=(1,1)$ is a rhombus point; 
see Fig. \ref{mos}. 

\begin{proof} 
By replacing in $\phi$ the points in $\mathcal{P}_1$,
a direct calculation shows that the equivalent $9 \times 8$ 
matrix has rank 7, where the null space of  $\phi$  is given by the vectors 
$(0, 0, 0, -2/3, 0, 0, 1, 0, 0)$ and
$(-2/3, 0, 0, 0, 1, 0, 0, 0, 0)$.
The linear combination of the  
corresponding polynomials leads to 
\begin{equation}
\label{polinomios-rombo}
f( {\tt a,d,}x,y)= 
a \left( 2x^3 - 3x^2 \right) 
+ 
d \left( 2y^3 - 3y^2 \right),
\ \ \ 
[a,\,d] \in \K \mathbb{P}^1.
\end{equation}
\end{proof}

\begin{remark} 
Behavior of the linear system at $\mathscr{A}$.
\begin{upshape} 
Let
$ \mathcal{P}=\{ (0,0), (1,0), (0,1), (x_4, y_4) \}$
be a configuration.

\noindent 
1.
If $(x_4, y_4)$ tends to be in a line 

\centerline{$
L_\alpha \subset \mathscr{A}
\big\{ \mathscr{A}\big\} \backslash 
\{ R_1=(1,1), R_2 =(-1,1), R_3= (1,-1)  \}, 
$}

\noindent 
then the polynomial $f(x_4, y_4,x,y)$ in 
\eqref{matriz-grado-tres}
has two lines of critical points in  
the respective pair of parallel $\K$--lines $L_\alpha$, $L_\beta$,
in 
the arrangement $\{\mathscr{A}(x,y) = 0\}$.
Figure \ref{atlas-polinomios-cubicos}
provides a sketch up to affine transformations.

\noindent
2. 
If $(x_4, y_4)$ tends to be
the vertex $(0,0) \in \triangle$,
then the polynomial $f(x_4, y_4,x,y)$ in \eqref{polinomio-grado-tres} becomes

\centerline{ 
$f(0,0, x,y)  
=
\frac{1}{3}(x^3 + y^3) - (x^2 y + xy^2) - 
\frac{1}{2} (x^2 + y^2) + xy$.
}

\noindent 
As is expected, the curve 
$
\{f(0,0, x, y)= 0 \} 
$
has a cusp of multiplicity two at (0,0),
see
Fig. \ref{atlas-polinomios-cubicos}.
The same is valid if $(x_4, y_4)$ tends to be
any other vertex
$(1,0)$, $(0,1)$ of $\triangle$. 
Figure \ref{atlas-polinomios-cubicos} 
shows $f(1,0,x,y)$, corresponding to $V_2=(0,1)$
denoted as ${\tt V}_2$ in the figure.
\end{upshape}
\end{remark}

\begin{remark}
\begin{upshape}
Let $\mathcal{P}$ be any 
configuration of four points. Thus 
$\mathcal{L}_3(\mathcal{P}) \neq \emptyset$, 
{\it i.e.}
there exists
a non constant degree three polynomial having critical 
points at least in $\mathcal{P}$. 
\end{upshape}
\end{remark}

\subsection{Affine classification of quadrilateral configurations}

We now study the independence of the previous results
\S4.1, with 
respect to the coordinate system.

A valuable tool in the study of  polynomials of degree 
three is the action of 
the group of affine automorphisms of $\K^2$,
say  $\aff(\K^2)$. 
It is a six $\K$--dimensional Lie group.
Let $\aff(\K^2)$ acts on the space of polynomials 
of degree $d$ as
\begin{equation}\label{accion-afin}
\aff(\K^2) \times \K[x,y]_{=d}  \longrightarrow 
\K[x,y]_{=d},
\ \ \ 
(T,f)\longmapsto f\circ T.
\end{equation}

\noindent 
This action is 
rich enough and yet treatable.
The affine group acts on configurations such as
\begin{equation}
\label{accion-afin-2}
\aff(\K^2) \times {\it Conf}(\K^2, n) 
\longrightarrow {\it Conf}(\K^2, n), 
\ \ \
(T,\mathcal{P}) \longmapsto  T^{-1}(\mathcal{P})
.
\end{equation}
Thus, 
if $f \in \K[x,y]_{=d}$ has $n$ isolated
critical points, say $\mathcal{P} \in {\it Conf}(\K^2, n)$, 
then 
$f \circ T$ has critical points at 
$T^{-1}(\mathcal{P})$. 
Hence, a useful associated 
object is the quotient space of
quadrilateral configurations up to affine transformations.

\begin{definition}
\begin{upshape}
The space of {\it generic quadrilateral configurations} is 
\begin{equation}
\label{cuadrilateros-genericos-1} 
\mathcal{Q}= 
\left\{ 
\mathcal{P}_0=
\{(x_{1\, 0}, y_{1 \, 0}), \ldots, 
(x_{4 \, 0}, y_{4\, 0})\}
\ \Big\vert \
\begin{array}{l}
\hbox{quadrilaterals configurations} 
\\
\hbox{having no three collinear vertices }
\\
\hbox{or determining two parallel lines}
\end{array} 
\right\} 
\subsetneq {\it Conf}(\K^2, 4).
\end{equation}
\end{upshape}
\end{definition}

Note that a quadrilateral configuration $\mathcal{P}_0$ 
does not have order. 
It determines several quadrilaterals, 
{\it i.e.}
with a cyclic order in its vertices. 
Let 
$$
\triangle= \{V_1=(0, 0), V_2 = (1,0), V_3=(0,1) \},
\ \ \
\Delta = \{ {\tt V}_1=(0, 0), {\tt V}_2 = (1,0), 
{\tt V}_3=(1/2, \sqrt{3}/2) \} 
$$

\noindent
be two triangles.
Consider a linear transformation  
$R \in GL(2, \K)$  such that  $R(\triangle) =\Delta$,
$R(V_2) = {\tt V}_2$ and $R(V_3) = {\tt V_3}$,  
see Fig. \ref{mos}. 
The affine symmetries of $\Delta$, 
\begin{equation}
Sym(3)=
\label{simetrias-de-Delta}
\{  
\sigma_\alpha \in  \aff(\K^2 ) 
\ \vert \  \sigma_\alpha(\Delta) = 
\Delta, \ \alpha \in 1, \ldots , 6 \},
\end{equation}

\noindent 
are isomorphic to the symmetric group of order 3;
with three reflections  
$\sigma_2, \sigma_4, \sigma_6$ (with axis 
in the lines ${\tt N}_1, \, {\tt N}_2, \, {\tt N}_3$)
and their products
$\sigma_1=id, \sigma_3, \sigma_5$;
see Fig. \ref{mos}.b.
By abusing the notation, $Sym(3)$ also denotes
the affine symmetries of $\triangle$.

\noindent 
Thus, we use three coordinate systems as follows. 
Let 
$ \mathcal{P}_0=
\{ (x_{1\,0}, y_{1\,0}), \ldots,  (x_{4\,0}, y_{4\,0})\} $
as in \eqref{cuadrilateros-genericos-1}.
By using the affine action, 
we reduce $\mathcal{P}_0$ to  
$\{(x_4, y_4) \}$ or   
$\{ ({\tt x}_4, {\tt y}_4) \}$.
There are affine maps 
$T_j \in \aff(\K^2)$,

\begin{equation}\label{normalizacion}
\begin{picture}(200,65)
\put(20,56){$
\mathcal{P}_0=
\{ (x_{1\,0}, y_{1\,0}), \ldots,  (x_{4\,0}, y_{4\,0})\} $}
\put(-64,4){$\mathcal{P}=\{ 
\underbrace{V_1, V_2, V_3}_{\triangle},  V_4= (x_4, y_4)\}$}
\put(128,4){$\{ 
\underbrace{{\tt V}_1, {\tt V}_2, {\tt V}_3}_{\Delta}, 
{\tt V}_4 = ({\tt x}_4, {\tt y}_4)\} \, .$}
\put(69,4){\vector(1,0){55}}
\put(75,49){\vector(-1,-1){33}}
\put(90,-8){$R$}
\put(88,11){$R^{-1}$}
\put(40,35){$T_j$}
\put(138,35){$R \circ T_j$}
\put(115,49){\vector(1,-1){33}}
\put(123,8){\vector(-1,0){55}}
\end{picture}
\end{equation}

\bigskip

\noindent
By notational simplicity, we also denote by $\mathcal{P}$
the configuration on the right side.

A key point is the number of affine maps $\{ T_j \}$,
depending on $\mathcal{P}_0$ 
to be computed in Corollary \ref{corolario-son-24}. 
 
\begin{figure}[h!]
\centerline
{\includegraphics[scale=0.35]{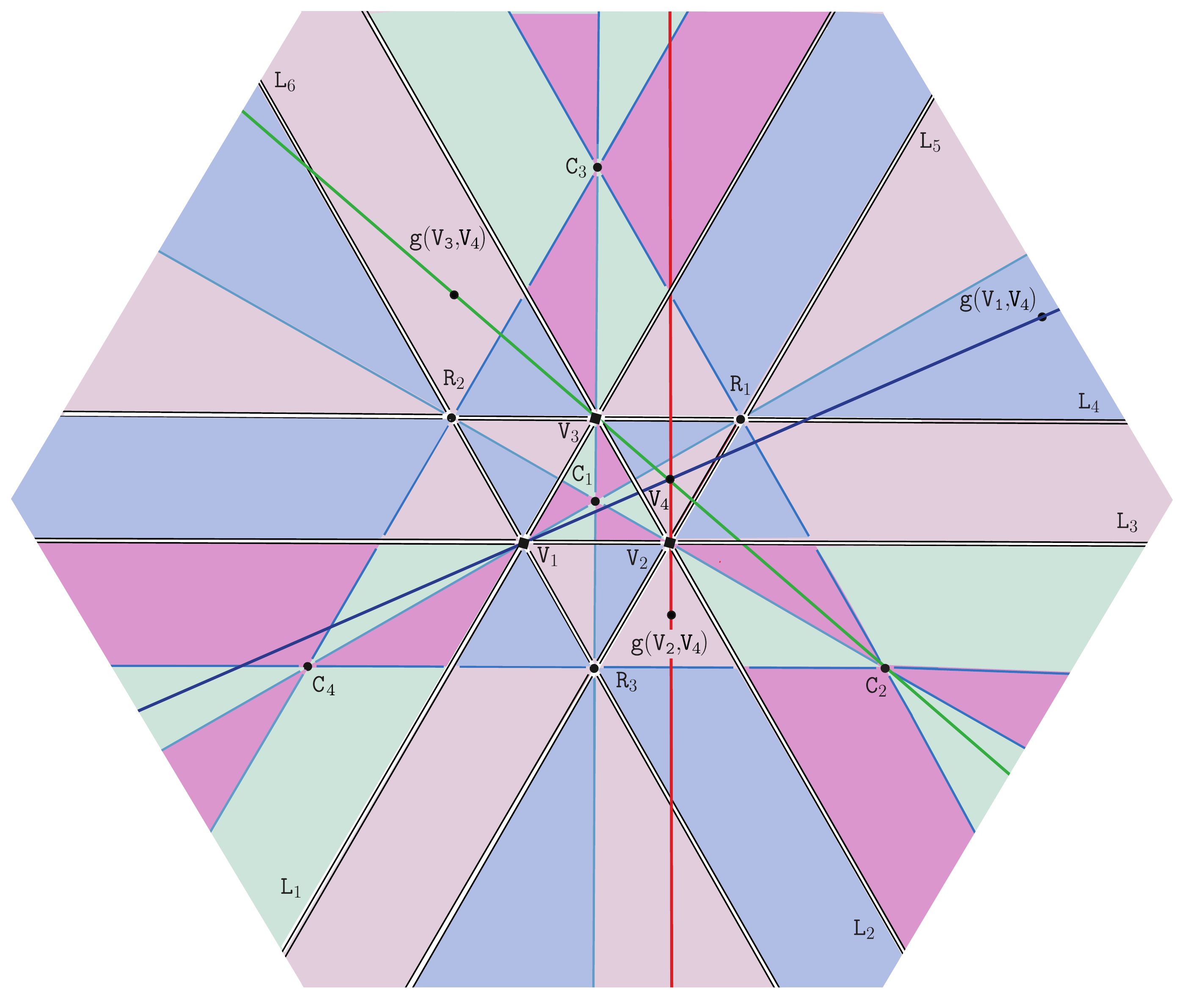}}
\caption{ 
The plane $\R^2 \backslash \tt A$ with 
coordinates $\{{\tt x}_4, {\tt y}_4 \}$ 
parametrizes the quadrilateral configurations
$\{ {\tt V}_2, {\tt V}_2,{\tt V}_3, {\tt V}_4=({\tt x}_4, {\tt y}_4) \}$.
The pair tile $\tt Q= {\tt Q}_1 \cup {\tt Q}_2$
is a fundamental domain for
the moduli space of quadrilateral configurations, up to $\aff(\K^2)$--equivalence. 
There are 24 copies of the fundamental region
${\tt Q}$.
We colored  ${\tt Q}_2$ and its copies 
pink or blue (resp. ${\tt Q}_1$ and its copies 
aquamarine or magenta)
tiles for strictly convex (resp. non convex) 
quadrilateral configurations. 
\label{transformaciones-gauge-2}
}
\end{figure}

In accordance with
Fig. \ref{mos} and  \ref{transformaciones-gauge-2}, 
the triangles $\triangle$, $\Delta$ determine
the points, line arrangements and regions below.

\smallskip 

\noindent
$\bigcdot$
Three {\it rhombus points} $R_1, \, R_2, \, R_3$
(resp. ${\tt R}_1, \ {\tt R}_2, \, {\tt R}_3$). 

\smallskip

\noindent
$\bigcdot$
Four {\it center points} $C_1, \, C_2, \, C_3, \, C_4$
(resp. ${\tt C}_1, \, {\tt C}_2, \, {\tt C}_3, \, 
{\tt C}_4$).

\smallskip

\noindent
$\bigcdot$
A six line arrangement 
$\mathscr{A}=L_1\cup \ldots \cup L_6$
(resp. ${\tt A}={\tt L}_1 \cup \ldots \cup {\tt L}_6$)
sketched as six double lines. $\mathscr{A}$ was
already described in the introduction and in 
\eqref{arreglo-de-seis-lineas}.

\smallskip

\noindent
$\bigcdot$
A six line arrangement 
$\mathscr{B}=N_1\cup \ldots \cup N_6$
(resp. ${\tt B}={\tt N}_1 \cup \ldots \cup {\tt N}_6$),
sketched as six blue lines, 
where $N_1, \, N_2, \, N_3$ are the axis of
symmetry of $\triangle$. 
The lines $N_1, \, N_2, \, N_3 $
are fixed under 
$\sigma_1, \, \sigma_2, \, \sigma_3$
in $\aff(\R^2)$ leaving invariant $\triangle$.
The lines 
$N_4, N_5, N_6$ determine the triangle $C_1, C_2, C_3$.  

\smallskip 

\noindent 
Naturally these points and arrangements are 
in correspondence under the map $R$ 
in \eqref{normalizacion}.

\smallskip
\noindent 
$\bigcdot$
In case $\K= \R$, 
we have two open connected regions in $\R^2$; 

\noindent 
convex  quadrilateral
configurations when $(x_4, y_4) \in Q_1$ (aquamarine) 
and 

\noindent non convex for 
$Q_2$ (magenta).

\noindent Analogously, we have
${\tt Q}_1=R(Q_1)$ and ${\tt Q}_2= R(Q_2)$.
Moreover, 
the boundary of $Q_1$, $Q_2$ shall be described by 
using the isotropy of the respective configurations.

\begin{lemma} Let $\mathcal{P} \in  \mathcal{Q}$
be a generic quadrilateral configuration in $\K^2$ as 
in \eqref{cuadrilateros-genericos-1}.
If the affine isotropy group of $\mathcal{P}$
$$
\aff(\K^2)_{\mathcal{P}}  
\doteq
\{ T \in \aff(\K^2) \ \vert \  
T^{-1}(\mathcal{P}) = \mathcal{P} \}
$$ 

\noindent 
is non trivial, then it is
isomorphic to one of the subgroups below.

\smallskip

\noindent 
Case 1. $\aff(\K^2)_{\mathcal{P}} \cong Sym(3)$ if and only if
up to affine transformation 
$\mathcal{P}$ has vertices in an 
equilateral triangle and its 
center. 

\smallskip

\noindent 
Case 2.  $\aff(\K^2)_{\mathcal{P}} \cong \ZZ_2 \times \ZZ_4$ if and only if
up to affine transformation $\mathcal{P}$ is a rhombus
(its vertices 
determine a pair of two parallel lines).

\smallskip

\noindent 
Case 3. 
$\aff(\K^2)_{\mathcal{P}} \cong \ZZ_2 $
if and only if 
up to affine transformation

\smallskip

\noindent 
i)  
$\mathcal{P} = \{(0,0), (1,0), (1/2, \sqrt{3}/2), 
({\tt x}_4,{\tt y}_4) \}$
where $ ({\tt x}_4, {\tt y}_4)$ is a fixed 
point under the reflection $\sigma^{\prime}_2$ with axis ${\tt N}_2$
in the isotropy  of the triangle $\Delta$
and it is different of the center of $\Delta$, or
 
\smallskip

\noindent 
ii) $\mathcal{P}$ is a trapezoid,
its vertices 
determine two parallel lines,  
different from a rhombus.
\hfill
$\Box$
\end{lemma}

\begin{corollary} 
\label{corolario-son-24}
Let $\mathcal{P}_0$ be a generic quadrilateral configuration, 
the following assertions are equivalent. 

\noindent  1)
$\mathcal{P}_0$ 
has a trivial isotropy group
$\aff(\K^2)_{\mathcal{P}_0} =id$. 

\noindent 
2) 
There are 24 affine transformations 
$R \circ T_j$ in 
\eqref{normalizacion}, sending $\mathcal{P}_0$
to 
$\{(0,0), (1,0), (1/2, \sqrt{3}/2), ({\tt x}_4, {\tt y}_4) \}$.
\hfill 
$\Box$
\end{corollary}

Now we compute 
the orbit $\{ R \circ T_j (\mathcal{P}_0) \}_{j=1}^{24}$
in terms of the fourth point in 
$\{({\tt x}_4, {\tt y}_4) \} \in \R^2$.
Certainly,
the orbit has obvious elements given by 
the affine symmetries of $\Delta$. The 
non intuitive transformations 
between quadrilateral configurations 
$R \circ T_j (\mathcal{P}_0)$, are computed in the following result.

\begin{lemma}\label{gauge-dependiendo-de-W}
Let 

\centerline{$\{ \underbrace{(0,0), (1,0), (1/2, \sqrt{3}/2),}_{\Delta} {\tt V}_4
=({\tt x}_4, {\tt y}_4) \}$}

\noindent
be a generic quadrilateral configuration and consider a vertex 
${\tt V}_j \in \Delta$.  
There exist three $\K$--rational diffeomorphisms 
(different from the identity)
\begin{equation}\label{gauge-equation}
{\tt g} ({\tt V}_{j}, \ \ ): \K^2 \backslash {\tt A}
\longrightarrow 
\K^2 \backslash {\tt A}
, \ \ \
{\tt V}_4 \longmapsto  {\tt g}({\tt V}_j, {\tt V}_4)
, 
\ \ \
j \in 1, 2, 3,
\end{equation}
such that 
the quadrilateral configurations 

\centerline{
$\{(0,0), (1,0), (1/2, \sqrt{3}/2), {\tt V}_4
\}$
\ \ and \ \ 
$\{(0,0), (1,0), (1/2, \sqrt{3}/2), 
{\tt g}({\tt V}_j, {\tt V}_4) 
\}$
}

\noindent 
are $\aff(\K^2)$--equivalent. 
\end{lemma}

We note that 
${\tt g}({\tt V}_j, \ \ )$ are non affine maps. 

\begin{proof}
The choice of one vertex ${\tt V}_j \in \Delta$, determines an opposite side $\Delta$. 
Without loss of generality, we consider the vertex 
${\tt V}_2= (1,0) \in \Delta$  and
${\tt L}_1 =\{{\tt y}- \sqrt{3}{\tt x}= 0\}  
\subset {\tt A}$ is the opposite side; 
see Fig. \ref{gauge}.

\noindent 
For fixed $j=2$, we consider ${\tt V}_4$.
Let ${\tt L}$ be the line by ${\tt V}_4$ and ${\tt V}_2$;
${\tt L}$ is the red line in Fig.  
\ref{gauge}.
We assume that ${\tt L}_1$ and ${\tt L}$ are non 
parallel. 
\noindent 
There exists a unique $\K$--affine embedding

\centerline{
$
\mathfrak{j} : \K \longrightarrow \K^2, 
\ \ \
\hbox{ with } \
\mathfrak{j}(\K)= {\tt L}, \ 
\mathfrak{j}(1) = {\tt V}_2, \
\mathfrak{j}(0) = {\tt L}_1 \cap {\tt L} 
\doteq {\tt 0} .
$}

\noindent
The definition of the map in ${\tt L}$
is 
\begin{equation}\label{gauge-puntual}
{\tt g}({\tt V}_2, \ \ ) :
{\tt L}\backslash \mathfrak{j}(0)
\longrightarrow 
{\tt L} \backslash \mathfrak{j}(0), 
\ \ \
{\tt V}_4 \longmapsto 
\mathfrak{j} 
\left( \frac{1}{ \mathfrak{j}^{-1} 
({\tt x}_4, {\tt y}_4)} \right) .
\end{equation}

\noindent 
Secondly, we shall extend this definition for 
${\tt V}_4 \in \K^2 \backslash {\tt L}_1$. 
In order to 
avoid cumbersome computations, the coordinates  
$\{ (x,y) \}$  in \eqref{normalizacion} are more
suitable. 
Assume $\mathcal{P}= \{(0,0), (1,0), (0,1), (x_4, y_4) \}$,
the vertex is $V_2=(1,0) \in \triangle$ and
$L_1 = \{x_4=0 \}$ is the opposite side.
The analogous definition provides the rational map  
\begin{equation}
\label{gauge-en-el-plano}
\begin{array}{rcl}
g(V_2, \ \ ):
\K^2 \backslash \{x_4 (x_4-1)= 0 \} & \longrightarrow  & 
\K^2 \backslash \{x_4 (x_4-1)= 0 \} , 
\\ 
V_4 = (x_4, y_4) & \longmapsto & 
\Big( \dfrac{1}{x_4}, \dfrac{-y_4 +y_4 x_4}{x_4 -1}\Big).
\end{array}
\end{equation}

\noindent
It enjoys the properties described below.

\noindent
$\bigcdot$ 
$g(V_2, \ \ )$ is a birational map of $\K^2$.

\noindent
$\bigcdot$ 
$g^{-1}(V_2, \ \ ) = g(V_2, \ \ )$,
{\it i.e.} it is an involution. 

\noindent 
$\bigcdot$
The point $V_2$ and the line $\{ x=-1\}$
are fixed under $g(V_2, \ \ )$.

\noindent 
$\bigcdot$
The poles of the map $g(V_2, \ \ )$ are localized at
$\{ x=0\}$ and 
$\{ x-1=0\} \backslash \{(0,1) \}$. Thus, strictly 
speaking the map 
is a $\K$--analytic diffeomorphism on 
$\K^2 \backslash \{ x(x-1)=0\}$. 
In the synthetic definition 
\eqref{gauge-puntual},
${\tt L}_1$ and ${\tt L}$ are non 
parallel. This originates the pole of $g(V_2, \ \ )$
at $\{x-1=0\}$. 

\noindent 
$\bigcdot$
A straightforward computations shows that
the line arrangements $\mathscr{A}$ and $\mathscr{B}$ 
(double and blue lines in Fig.
\ref{transformaciones-gauge-2}) 
are poles or remain invariants under 
$g_2(V_2, \ \ )$. 

Summarizing, we define 
\eqref{gauge-puntual} as

\centerline{$
{\tt g}({\tt V_2}, \ \  ) \doteq 
R \circ g(V_2, \ \ ) \circ R^{-1} 
$.}

Finally, given ${\tt V}_4$ and 
${\tt g}({\tt V_2}, {\tt V}_4)$,
there exists a unique transformation 
$T \in \aff (\K^2)$,
which leaves the line ${\tt L}_1$ fixed 
so that
$T({\tt V}_4) = {\tt g}({\tt V_2}, {\tt V}_4)$;
see Fig. \ref{transformaciones-gauge-2}.
Under $T$, the quadrilateral configurations
$$
\{ (0,0),  (1,0), (1/2, \sqrt{3}/2), {\tt V}_4 \}
\hbox{ \ \ and \ \ } 
\{ (0,0), (1,0), (1/2, \sqrt{3}/2), T ({\tt V}_4) \}
$$

\noindent  
are affine equivalent.

\begin{figure}[h!]
\centerline
{\includegraphics[scale=0.40]{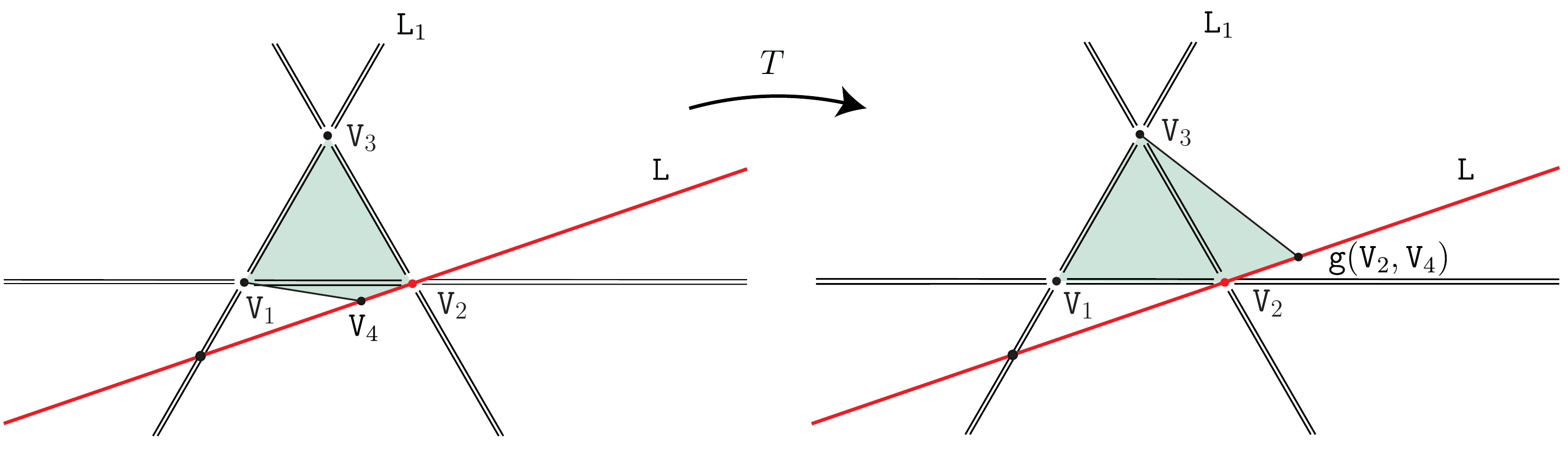}}
\caption{ 
The point ${\tt g} ( {\tt V}_2, {\tt V}_4 )$ determines
an affine map $T$ between generic quadrilateral configurations.
\label{gauge}
}
\end{figure}

The other vertices of 
the triangle $\Delta$ 
determine rational maps 
$ {\tt g}({\tt V}_1, \ \ ), \,  
{\tt g}({\tt V}_3, \ \ )$, both
enjoy analogous properties. 
\end{proof}

\begin{remark}
\begin{upshape}
Three blue lines in Fig. \ref{transformaciones-gauge-2}
correspond to the fixed points under the 
reflection symmetries
$Sym(3)$ of $\Delta$. By using \eqref{gauge-puntual}, the complete 
configuration of six blue lines $N_1, \ldots, N_6$ is invariant under
the three transformations ${\tt g}({\tt V}_j, \ \ )$.
We leave this 
assertion for the reader. 
\end{upshape}
\end{remark}

\begin{lemma}\label{moduli-cuadrilateros} 
1. The quotient space of generic quadrilateral configurations up to affine transformations, given by 
\begin{equation} 
\label{proyeccion-moduli}
\pi: \mathcal{Q} \longrightarrow 
\mathcal{Q} / \aff(\K^2)
, \ \ \
\{(x_{1\,0}, y_{1\, 0}), \ldots , (x_{4\, 0}, y_{4\, 0} ) \}
\longmapsto [ ({\tt x}_4, {\tt y}_4 ) ],
\end{equation}
is a $\K$--analytic surface $\tt Q$.

\noindent 2. 
For $\K = \C$, the quotient $\tt Q$ is a connected complex surface.

\noindent 3. 
For $\K = \R$, the quotient  
has two connected components
$\tt Q = {\tt Q}_1 \cup {\tt Q}_2$
and singular points with local 
models $\K^2 / \mathbb{Z}_2$ 
or $\K^2 /Sym(3)$.  
\end{lemma}

Some comments are in order.  
Figure \ref{transformaciones-gauge-2} illustrates 
the fundamental domains for $\pi$ over $\K=\R$. 
The double lines 
${\tt A}= {\tt L}_1 \cup \ldots \cup {\tt L}_6$ 
in Fig. \ref{mos}--\ref{atlas-polinomios-cubicos} correspond to
forbidden positions for $({\tt x}_4 ,{\tt y}_4)$.
Moreover, 
$({\tt x}_4 ,{\tt y}_4) \in {\tt Q}_1$ means 
a non convex quadrilateral configuration;
$({\tt x}_4 ,{\tt y}_4) \in {\tt Q}_2$ determines a 
strictly convex quadrilateral configuration. 

\begin{proof}
The set theoretical construction of the quotient is simple,
and we describe its projection 
$\pi$ in \eqref{proyeccion-moduli}.
Given $\mathcal{P}_0 \in \mathcal{Q}$,
we apply an affine transformation $R \circ T_j$  in 
\eqref{normalizacion}
sending it to

\centerline{
$ R \circ T_j(\mathcal{P})=
\{(0, 0), (1, 0), (1/2,  \sqrt{3}/ 2), 
{\tt V}_4=({\tt x}_4, {\tt y}_4) \} 
$.}

\noindent 
Case 1. The isotropy is trivial 
$\aff(\K^2)_{\mathcal{P}} =id$. 
There are exactly $24$ different choices for 
$R\circ T_j$, as in Lemma \ref{corolario-son-24};
we have that $\pi$ has as target 
$\K^2=\{({\tt x}_4,  {\tt y}_4)\}$.  

In order to describe its analytic properties, 
recall that the Klein four--group $K$ is isomorphic to 
$\mathbb{Z}_2 \times \mathbb{Z}_2$. 
It is such that each element is self--inverse 
(composing it with itself produces the identity) 
and  composing any two of the three 
non--identity elements produces the third one;
see \cite{Amstrong} p.\,87. 
Moreover, the group $Sym(4)$ is of order 24,  
having a Klein four--group $K$ as a proper normal subgroup; 
thus $Sym(3) = Sym(4)/K$.  
We recognize 
$$
K =\{ id, \,  {\tt g}({\tt V}_j, \ \ ) \ \vert \ j \in 1,2,3  \}
$$
as the group in Lemma \ref{gauge-dependiendo-de-W}. 
Recall \eqref{simetrias-de-Delta} and consider the homomorphism 
given by
$$
\varphi: Sym(3) \longrightarrow Aut(K), 
\ \ \
\sigma \longmapsto   
\sigma^{-1}_\alpha \circ  {\tt g}({\tt V}_j, \ \ ) \circ  
\sigma_\alpha  ({\tt x}_4, {\tt y}_4).
$$
The semidirect product of $K$ and $Sym(3)$ determined by 
$\varphi$ is 
$
Sym(4) = K \rtimes_\varphi Sym(3), 
$
see \cite{Amstrong} p.\,133.
Hence we have a representation of $Sym(4)$ in the 
birational 
transformations of $\K^2\backslash {\tt A}$ and 
\begin{equation}\label{cociente}
\tt Q = \frac{\mathcal{Q}}{\aff(\K^2)} = 
\frac{\K^2 \backslash {\tt A} }{Sym(4)}
\end{equation}

\noindent 
is the quotient space. 
See \cite{Prill} for a general theory of the quotients 
of complex manifolds under a discontinuous group of 
automorphisms. Assertion (1) is done.

For assertion (2), $\K=\C$;
note that $\K^2 \backslash {\tt A}$ is a connected complex manifold. 
The local behavior of this complex quotient at the points with non trivial 
isotropy  $\ZZ_2$ at the lines 
${\tt N}_1, \, {\tt N}_2, \, {\tt N}_3$ is known to 
be non--singular
(because of C. Chevalley \cite{Chevalley}, see also \cite{Flatto}).
For ${\tt C}$ the isotropy is $Sym(3)$ and the same references describe
the local structure of the quotient.

For assertion (3), $\K=\R$;  
clearly the convexity or non convexity 
of a quadrilateral configurations
are affine invariants, whence there are two connected components. 
At the points 
${\tt C}, \ldots ,{\tt C}_4$ and lines 
${\tt N}_1, \, {\tt N}_2, \, {\tt N}_3$
where 
the isotropy of the quadrilateral configurations 
is non trivial, the quotient 
\eqref{cociente}
has singularities; it is an orbifold.  
\end{proof}

As final step in the proof of Theorem \ref{segundo-teorema-grado-tres},
we consider the action on projective classes
\begin{equation}\label{LaAccion}
\mathcal{A}:\aff(\K^2) \times Proj(\K[x,y]_{=3}) \longrightarrow 
Proj(\K[x,y]_{=3}),
\ \ \ (T, [f]) \longmapsto [f \circ T].
\end{equation}

\noindent
This action provides
an $\aff(\K^2)$--bundle structure on $\K[x,y]_{=3}$. 
Denote the 
\emph{stabilizer} or \emph{isotropy group} of $[f] \in Proj(\K[x,y]_{=3})$ by 
$$
\aff(\K^2)_{[f]} \doteq \{T\in \aff(\K^2) 
\ \vert \ 
f \circ T = \lambda f, \ \lambda \in \K^* \}.
$$
\noindent 
Equations
\eqref{solucion-general-grado-3-usando-base-canonica}
and \eqref{normalizacion}
provide bijective correspondence between 
generic quadrilateral configuration 
in $({\tt x_4}, {\tt  y_4}) \in \K^2 \backslash {\tt A} $
and projective classes of polynomials
$[f( R^{-1}({\tt x}_4  ,{\tt y}_4),x,y)]$.
If $\mathcal{P} \in {\tt Q}$, then we verify that 
the isotropy of the quadrilateral configuration
$\aff(\K^2)_\mathcal{P}$ is isomorphic to
$\aff(\K^2)_{[f]}$.  
Thus, we have a section
$$
f \circ R^{-1}: \K^2\backslash \{ {\tt A} \} 
\longrightarrow Proj(\K[x, y]_{=3}), 
\ \ \ 
({\tt x}_4, {\tt y}_4) \longmapsto 
[f( R^{-1}({\tt x}_4  ,{\tt y}_4),x,y)]
$$
\noindent 
and a diagram 
\begin{equation}\label{bundle-grado-tres}
\begin{picture}(200,55)
\put(133,46){$Proj(\K[x, y]_{=3})$}
\put(144,39){\vector(0,-1){33}}
\put(148,21){$\pi$}
\put(-16,20){$[f( R^{-1}({\tt x}_4  ,{\tt y}_4 ),x,y)]$}
\put(20,-13){$\K^2\backslash  \{{\tt A}  \}$}
\put(45,0){\vector(2,1){80}}
\put(120,-13){$\dfrac{Proj(\K[x, y]_{=3})}{\aff (\K^2)}$\,,}
\end{picture}
\end{equation} 
\medskip

\bigskip

\noindent 
here $\pi$ is the projection of classes from the action \eqref{LaAccion}. The $\aff(\K)$--orbit of a projective class 
$[f] \in \K[x,y]_{=3}$ is homeomorphic
to $\aff(\K^2) / \aff(\K^2)_{[f]}$. 
Obviously,
$\K[x,y]_{=3, id}$ is open and dense in $\K[x,y]_{=3}$.

The proof of assertion 1, Theorem 1 is done.
 
\begin{remark}
\begin{upshape}
It is well known (see for instance 
\cite{Duistermaat-Kolk} p.~53), 
that if we consider 
$$
\K[x,y]_{=3,id} \doteq \{  f \in \K[x,y]_{=3} \ \vert \  \aff(\K^2)_{f} ={id}\}
$$
then 
the restricted action in 
$\K[x,y]_{id} $, determines
a principal fiber $\aff(\K^2)$--bunldle structure.
In particular, 
the quotient $\K[x,y]_{=3, id}/\aff(\K^2)$ 
is a two dimensional $\K$--analytic manifold.
\end{upshape}
\end{remark}

\begin{remark}
\begin{upshape}
For $\K= \R$, the fundamental domain 
${\tt Q}_1 \cup {\tt Q}_2$ determines the
bifurcation diagram of the respective 
Hamiltonian vector fields, 
see Fig. \ref{atlas-polinomios-cubicos}.
By construction, 
${\tt Q}_1 $ 
has two boundaries and one vertex ${\tt C}$
and
${\tt Q}_2$ has one boundary (without extreme points).

\begin{figure}[h!]
\centerline
{\includegraphics[scale=0.5]{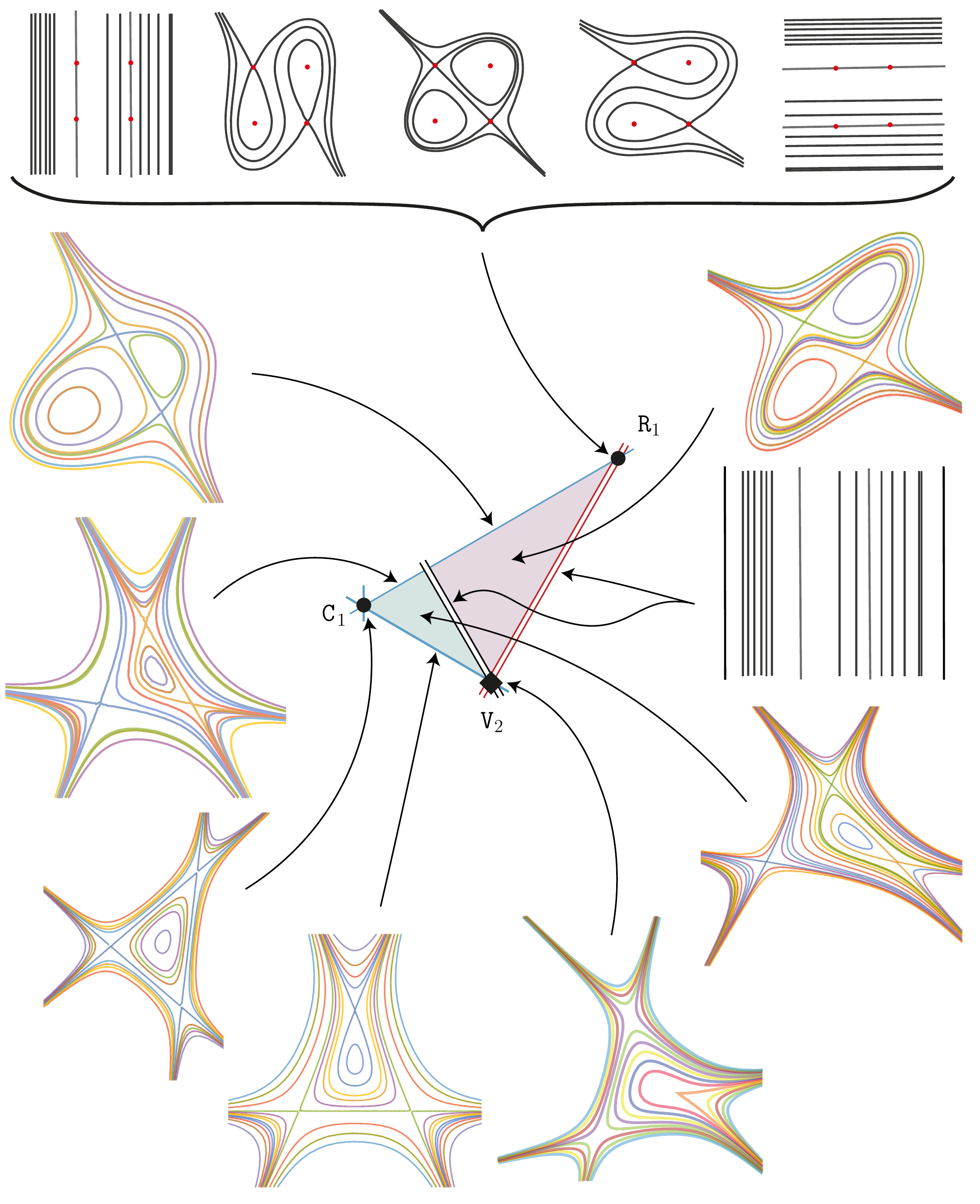}}
\caption{
Bifurcation diagram of the real Hamiltonian vector fields
$X_{f \circ R^{-1}}$ according to  the position of  
four critical points in the fundamental region $\tt Q$. 
At the rhombus point $\tt R_1$,
the configuration of four points
$\mathcal{P}=\{(0,0), (1,0), (0,1), {\tt R}_1=(1,1) \} \subset \Sigma(f_\theta)$ is common; 
see Example \ref{familia-de-rotaciones-texto}.
The upper row illustrates
the topology of $\{ f_\theta (x,y) \ \vert \ \theta \in [0,\pi /2 ]\}$.  
A saddle connection bifurcation occurs for
$\theta= \pi / 4$. 
See \url{https://github.com/alexander-arredondo/Mathematica-code-for-Essentially-determined-polynomials-of-degree-3/commit/e6a08f9a20da7b23d7a72beff8290af3a23260dc} for a code animation in Mathematica of this situation.
\label{atlas-polinomios-cubicos}
}
\end{figure}
\end{upshape}
\end{remark}


We summarize the results in Table 2.

\begin{table}[htbp]\label{dimensiones-para-4-puntos}
\begin{center}
\begin{tabular}{|c|c|c|c|c|}
\hline
&&&&$\vspace{-.4cm}$
\\
configuration $\mathcal{P}$ 
& cardinality & 
$dim_\K (\mathcal{L}_{3}(\mathcal{P}))$
& generators  & isotropy
\\
$\{(0,0), (1,0), (0,1), (x_4, y_4) \} $
& of $\Sigma(f) \vspace{-.3cm}$ &  
& of  $\mathcal{L}_{ 3}(\mathcal{P})$ & 
$\aff(\K^2)_{\mathcal{P}}$
\\
&&&&$\vspace{-.1cm}$
\\
 \hline
&&&&$\vspace{-.3cm}$
\\ 
$(x_4, y_4)\in Q$& 4 & 0 & 
eq.  \eqref{solucion-general-grado-3-usando-base-canonica-0} & $id$
\\
&&&& $\vspace{-.3cm}$ \\
\hline
&&&&$\vspace{-.3cm}$\\
$(x_4, y_4)=(1/3, 1/3) = C_1 $ &4&0& $xy(y+x-1)$ & $Sym(3)$
\\
&&&& $\vspace{-.3cm}$ \\
\hline
&&&&$\vspace{-.3cm}$\\
$(x_4, x_4), \,  x_4 \neq 0, \, 1$& 4 
& 0 & eq. \eqref{solucion-general-grado-3-usando-base-canonica-0} & $\mathbb{Z}_2$
\\
&&&&$\vspace{-.3cm}$
\\
\hline
&&&& $\vspace{-.3cm}$\\
$(x_4, y_4)=(1,1)= R_1$  & 4 & 1 $\vspace{-.03cm}$ & 
$\begin{array}{c}
2y^3-3y^2, \ 2x^3-3x^2
\\
\hbox{eq. \eqref{polinomios-rombo}}
\end{array}$  
& $ \mathbb{Z}_2 \times \mathbb{Z}_4$
\\
&&&&$\vspace{-.3cm}$
\\
\hline
&&&&$\vspace{-.3cm}$
\\
$(1, y_4 ) \in L_5, \ y_4 \neq \pm 1$& $\infty$ 
& 0 & 
$2x^3 - 3 x^2$  
& $\aff(\K)$ 
\\
&&&& $\vspace{-.3cm}$\\
\hline
&&&&$\vspace{-.3cm}$\\
$\{ (0,0),(1,0),(x_3,0),(x_4,0) \}$  & $\infty$ & 
2 $\vspace{-.3cm}$ & $
\begin{array}{c}
y^3, xy^2, y^2
\\
\hbox{ Lemma \ref{casos-degenerados}.2}
\end{array}$   
 & $\mathbb{Z}_2$
\\
&&&&\\
\hline
&&&&$\vspace{-.3cm}$
\\
$\{ (0,0),(1,0),(0,1), (0,0) \}$
& 3,  4  or  $\infty$ & 2  & 
$
\begin{array}{c}
x^3 - 3 x^2, \ y^3 - 3 y^2, 
\\
x^2y+xy^2-xy 
\end{array}
$
& $\mathbb{Z}_2$
\\
\hline 
\end{tabular}
\vskip0.2cm 
\caption{Dimension, 
generators and isotropy for 
$\mathcal{L}_3(\mathcal{P})$, 
where $\mathcal{P}$ is a configuration with 4 points
(3 simple points and a double one in the last row).
}
\end{center}
\end{table}

\begin{example}
Relation to the classification of plane cubic curves. 
\begin{upshape}
The Hesse pencil of cubic curves is 
$$
\{ z^3 + x^3 + y^3 - 3 \mu zxy = 0\}, 
\ \hbox{resp. } \
\{ x^3 + y^3 - 3 \mu xy  +1 = 0\}, 
\ \ \ 
\mu \in \C^*, 
$$
in the projective plane
$\C \mathbb{P}^2= \{ [z,x, y] \}$, 
resp. the affine plane; 
see \cite{Artebani-Dolgachev}. 
The key property is that
any non singular plane cubic is projectively 
equivalent to a member of the Hesse pencil.  
The critical points of the affine Hesse polynomial

\centerline{$f(\mu, x, y )= x^3 + y^3 - 3 \mu xy  +1$}

\noindent 
determine a generic quadrilateral configuration  

\centerline{
$\left\lbrace 
(0,0), 
(\mu,\mu), 
\big(-\zeta_1 \mu,\, \zeta_2 \mu\big), 
\big(\zeta_2 \mu, \, -\zeta_1 \mu\big) 
\right\rbrace 
\subset 
\C^2 \backslash \R^2, 
$}

\noindent 
where $\{ 1, \zeta_2, \zeta_3 \}$ are the cube roots of unity.
In order to translate it to our language, up to the linear transformation 
$
M_\mu: \C^2  \longrightarrow  \C^2,
\ 
(x,y) \longmapsto   
\left( \mu x - \zeta_2 \mu y, 
\mu x + \zeta_3 \mu y \right) $.
The
quadrilateral configuration changes to

\centerline{$
\mathcal{P}=
\{ (0,0), (1,0), (0,1), 
( 2 \zeta_1 \mu^2, 
\big(1 + \zeta_2 \big)\mu^2) \}. 
$}

\noindent  
By Theorem \ref{segundo-teorema-grado-tres},
the affine Hesse polynomial

\centerline{$
f(\mu, \ \ , \ \ ) \circ M (x,y) =
\mu^3 \left( 2 x^3 - 3 x (-1 + y) y - 3 x^2 (1 + y) + y^2 (-3 + 2 y) \right)+1
$}

\noindent 
is essentially determined. 
Since these quadrilateral configurations are non real,  
they are different from 
the given in Fig. \ref{atlas-polinomios-cubicos}. 
\end{upshape}
\end{example}

\subsection{Non essential determined polynomials of degree three}

By completeness, we describe the polynomials
arising from the configurations
 
\centerline{$\mathcal{P}= \{(0,0), (1,0), (0,1), (x_4, y_4) \} 
\in {\it Conf}(\K^2, 4)$, \
$(x_4, y_4) \in \mathscr{A}$.  
}

\begin{lemma}\label{casos-degenerados}
1. Let $\mathcal{P}=
\{ (0, 0), (1, 0), (x_3, 0), (x_4, y_4) \}$, with
$x_3 \neq 0, 1$ and
$y_4\neq 0$,
then $dim_\K(Proj(\mathcal{L}_3(\mathcal{P})))= 0$.

\smallskip

\noindent 2.
Let $\mathcal{P}=
\{ (0, 0), (1, 0), (x_3, 0), (x_4, 0) \}$  
be a configuration
then $dim_\K(Proj(\mathcal{L}_3(\mathcal{P})))= 2$.
\end{lemma}

\begin{proof}
In assertion (1), 
up to an affine transformation 
we can assume $y_4=1$.
The corresponding 
cubic polynomial takes the form 
$
f(x,y)= a_4 
\left( 2y^3-3y^2 \right)$, 
where  $a_4 \in \K^*$.

\noindent 
For assertion 2,
we search for 
polynomials $f(x,y) \in \K[x,y]^0_{\leq 3}$  
with at least four  affine  collinear critical points. 
The matrix of Eq.
(\ref{matriz-grado-tres}) results in
the cubic polynomials 

\centerline{
$
f(x,y)= a_3 x y^2+a_4y^3+a_7y^2 = y^2(a_3 x+a_4 y + a_7),
\ \ \
[a_3, a_4, a_7] \in \K \mathbb{P}^2, 
$}

\noindent 
with a line of critical points in $\{ y=0\}$.
\end{proof} 

\begin{example}
\begin{upshape}
The elementary methods provide an insight in the case 
of a double point in $\Sigma(f)$. Let
$\mathcal{P}_2=\{ (0,0),(1,0),(0,1), (0,0) \}$ be such a configuration.
A basis for $\mathcal{L}_3 (\mathcal{P}_2)$ is

\centerline{$
x^3 - 3 x^2, \ y^3 - 3 y^2, 
\ x^2y+xy^2-xy. 
$}

\noindent 
The first and second polynomials have lines of 
singularities, 
the third one four isolated critical points. 
The family of polynomials is 
$$
f(a_1,a_2, a_4, x,y )=
a_1 (x^3 - 3 x^2)+
a_2 ( x^2y+xy^2-xy  )+
a_4(y^3 - 3 y^2),
\ \ \ 
[a_1, a_2, a_4] \in \K \mathbb{P}^2.
$$
As is expected, for values 
$\left\{ (a_1, a_2, a_4=a_2^2 /9 a_1 \right\}$
the two dimensional family
$f(a_1,a_2, a_4, x,y)$ determines polynomials 
with three isolated singular points, 
one of them of multiplicity two, 
see Fig. \ref{atlas-polinomios-cubicos}.
\end{upshape}
\end{example}


\section{Degree four polynomials }
\label{polinomios-grado-cuatro}

Let 
\begin{equation}\label{cuartica-general}
f(x,y)=a_1 x^4+a_2 x^3 y+ \ldots + a_{13} x+a_{14} y
\, 
\in \K[x,y]^0_{\leq 4}
\end{equation}

\noindent 
be a polynomial as in \eqref{polinomio-general-grado-d}.
Here 
by notational simplicity we have avoided the double subindex,
and  
let $\mathcal{P} = \{ (x_\iota, y_\iota) \ \vert  \ \iota \in 1, \ldots, 7  \}$
be a configuration of seven points. 
The associated linear system for Eq. \eqref{cuartica-general} is
\begin{equation}\label{matriz-grado-cuatro}
\left(
\begin{smallmatrix}
& & & & &\vdots& & & & &  & & & \\
 4x_\iota^3 & 3x_\iota^2y_\iota^2 & 2x_\iota y_\iota^2 & y_\iota^3 & 0 & 3x_\iota^3 & 2x_\iota y_\iota & y_\iota^2 & 0 & 2x_\iota & y_\iota & 0 & 1 & 0 \\
0 & x_\iota^3 & 2x_\iota^2y_\iota & 3x_\iota y_\iota^2 & 4y_\iota^3 & 0 & x_\iota^2 & 2x_\iota y_\iota & 3y_\iota^2 & 0 & x_\iota & 2y_\iota & 0 & 1 
\\
&&&&&\vdots&&&&&&&&
\end{smallmatrix}
\right) 
\left(
\begin{smallmatrix}
a_1 \\
\vdots \\
a_{14}
\end{smallmatrix}
\right)
=  \overline{0}, \ \ \ \iota = 1, \ldots, 7.
\end{equation}

\noindent 
The interpolation matrix $\phi$, Eq. 
\eqref{matriz-grado-cuatro}, is square.
Hence, for an open and dense set of configurations
$\{ \mathcal{P}\} \subset {\it Conf}(\K^2, 7)$
such that $\{ det(\phi)=0 \}$, 
the resulting
space of polynomials of degree four 
having these $\mathcal{P}$ as critical points
is empty.
In order to overcome this situation,  
we introduce the following concept.

\begin{definition} 
\begin{upshape}
Let $\K[x,y]^0_{\leq d}$ with an even dimension 
and 
$\delta(d) = \frac{1}{4} \left( d^2 + 3d  \right)$
as in \eqref{valores-de-delta-explicitos}. 
Given a configuration  
$\mathcal{P}_0 \in {\it Conf}(\K^2, \delta(d)-1)$, 
consider a point 
$(x, y) \in \K^2$ and 
$$
\mathcal{P}_1=
\Big\{
\underbrace{(x_1, y_1), \ldots , (x_{\delta(d)-1}, 
y_{\delta(d)-1}) }_{
\mathcal{P}_0 }, 
(x, y) 
\Big\}\, 
\in {\it Conf}(\K^2, \delta(d)).
$$
\noindent 
The {\it interpolation algebraic curve 
of $\mathcal{P}_0$}  is  
$$
\mathcal{I} =
\left\{ det \big( {\phi} ( x_1, y_1, \ldots, 
x_{\delta(d)-1}, y_{\delta(d)-1}, x, y) \big)= 0 \right\}
\ \hbox{  in } \K^2.
$$
\end{upshape} 
\end{definition}

Obviously, $\mathcal{I}$ depends on 
$\mathcal{P}_0$, 
by notational simplicity we omit this dependence. 
Thus, we have a map
\centerline{
$
\mathcal{P}_0 = 
\{ (x_1, y_1), \ldots , (x_{\delta(d)-1}, y_{\delta(d)-1} ) \}
\longmapsto \mathcal{I}. 
$}

\begin{proposition}
\label{curvas-de-interpolacion}
Let $\K[x,y]^0_{\leq d}$ having even dimension.

\smallskip

\noindent 
1. The interpolation curve $\mathcal{I}$ of $\mathcal{P}_0$
describes the position of the $\delta(d)$--th point 
such that 
$dim_\K (\mathcal{L}_{d} (\mathcal{P}_1)) \geq 0 $.

\smallskip 

\noindent 
2.  
There exists a Zariski open set 
$\{ \mathcal{P}_0 \} \subset {\it Conf}(\K^2, \delta(d)-1)$
such that the associated  $\{ \mathcal{I} \}$
are algebraic curves of degree $2d-2$ in $\K^2$.
\end{proposition}

\begin{proof}
For assertion (2),  we
consider the degree $d$ polynomial

\centerline{$
f(x,y)=a_1 x^d+a_2 x^{d-1} y+ \ldots + 
a_{\delta(d)-1} x+a_{\delta(d)} y.
$}

\noindent 
After fixing the configuration $\mathcal{P}_0$, the associated linear system only has free variables $x$, $y$,
and 
the linear system is as follows
\begin{equation}
\left(
\begin{smallmatrix}
& & & & &  \vdots& & & & &  & & & \\
 (d)x^{d-1} & (d-1)x^{d-2}y & (d-2)x^{d-3} y^2 & \cdots & 0 & (d-1)x^{d-2} & \cdots & y^2 & 0 & 2x & y & 0 & 1 & 0 \\
0 & x^{d-1} & 2x^{d-2}y & \cdots & 4y^3 & 0 & x^2 & 2x y & 3y^2 & 0 & x & 2y & 0 & 1 
\\
&&&&&\vdots&&&&&&&&
\end{smallmatrix}
\right) 
\left(
\begin{smallmatrix}
a_1 \\
\vdots \\
a_{\delta (d)}
\end{smallmatrix}
\right) = \overline{0}.
\end{equation}

\noindent 
The determinant of this matrix has  
$x^{2d-2}$ as higher degree monomial, we are done.
\end{proof}

We describe some interpolation curves $\mathcal{I}$. 

\begin{example}
\begin{upshape}
Let $f \in \K[x,y]^0_{\leq 4}$ be a polynomial 
having degree four and let $\mathcal{P}_0 = 
\{ (x_\iota, y_\iota) \ \vert  \ \iota \in 1, \ldots,6 \}$ 
be a fixed configuration of six different critical points of $f$. 

\noindent 
1. 
If  three points of  $\mathcal{P}_0$ are in a line 
$\{ x= 0 \}$
and  two points are in $\{x=1 \}$,  
then the interpolation curve $\mathcal{I}$, 
of $\mathcal{P}_0$, is given by
\begin{equation}\label{interpolacion-3paralelas}
\mathcal{I}(x, y) =\left( -1152 y_4^2 y_5^2 (y_4-1)^2 x_6 (x_6-1)\right) x (x-1) (x - x_6) g(x, y). 
\end{equation}

\noindent 
The 
$\mathcal{I}$ is reducible and singular, 
it is the product of three parallels lines  
and a  polynomial $g(x, y)$ that pass through  
the six points in $\mathcal{P}_0$.

\noindent 2. 
Let
$\mathcal{P}_0 = \{ (x_\iota, y_\iota) \ \vert  \ \iota \in 1, \ldots,6 \}$ 
be any configuration of six points 
in the grid of nine points 

\centerline{
$\mathcal{G}=\{x(x-1)(x-c_1)=0\} \cap \{y(y-1)(y-c_2)=0 \}$,
\  where 
$c_1, \, c_2 \notin \{0,1\}$.}

\noindent  
Therefore, 
the interpolation curve $\mathcal{I}$, associated to the  
seventh point $(x_7, y_7)$, is the product of
the six lines defining 
$\mathcal{G}$.

\noindent 3.
Let $\mathcal{P} = 
\{ (x_\iota, y_\iota) \ \vert  \ \iota \in 1, \ldots,6 \}$ 
be a configuration of six critical points of $f$. 
If the six points are distributed in a conic $Q$, 
then the interpolation curve $\mathcal{I}$, 
associated to the  seventh point $(x_7, y_7)$, 
contains the conic. That is, 
$\mathcal{I}=Qg$ 
for some $g \in \K[x,y]^0_{\leq 4}$.
\end{upshape}
\end{example}

A complete study of the interpolation curves $\mathcal{I}$
arising from configurations of six points is the goal of a future project.


\section{Pencils of polynomial vector fields}

Now we will consider some special configurations 
of $(d-1)^2 \geq 4$ points.

\begin{definition}
\begin{upshape}
Let $\{ F(x,y)=0 \}, \ \{ G(x,y)=0 \}$ be two algebraic
curves in $\K^2$, both of degree $d-1$ $(\geq 2)$.
We assume that they have  transversal intersection  in 
exactly $(d-1)^2$  affine points, thus 
\begin{equation}\label{puntos-base-de-pincel}
\mathcal{P}_{ci}=\{F(x,y)=0\} \cap \{ G(x,y)=0\} \in 
{\it Conf}(\K^2, (d-1)^2) 
\end{equation}
is a {\it complete intersection configuration}.
The associated {\it pencil of curves} is
\begin{equation}\label{pincel-curvas}
\left\{ 
\mu F(x,y) + \nu G(x,y) =0   \ \vert \ [\mu, \nu] \in \K\mathbb{P}^1   
\right\}.
\end{equation}

\noindent 
$\mathcal{P}_{ci}$
is the {\it base locus} of the pencil of curves. 
\end{upshape}
\end{definition}

Moreover, 
the choice of an ordered pair of polynomial functions
from \eqref{pincel-curvas}, not just curves say

\centerline{
$
\big( 
{\tt a}  F(x,y) + 
{\tt b}  
G(x,y), \ 
{\tt c} F(x,y) +  
{\tt d}  
G(x,y) \big)
$}

\noindent 
determines a {\it $SL(2, \K)$--pencil of polynomial 
vector fields}
\begin{equation}
\label{pincel-de-campos}
\mathfrak{F}(\mathcal{P}_{ci})=
\left\{ \
X_{\tt M}=
-\big( 
{\tt c}  
F(x,y) + 
{\tt d}  
G(x,y) \big) \del{}{x} + 
\big( 
{\tt a} 
F(x,y) +  
{\tt b}  
G(x,y) \big) \del{}{y}    
\ \Big\vert \
{\tt M}=
\left(
\begin{smallmatrix}
\tt{-c} & \tt{ -d} \\
\ \tt{a} & \ \tt{b} 
\end{smallmatrix}
\right) \in SL(2, \K )
 \ \right\}.
\end{equation}

\noindent 
The condition $det({\tt M}) \neq 0$ 
is equivalent with the fact 
the zero locus
of the vector field $X_{\tt M}$
coincides with $\mathcal{P}_{ci}$.

\begin{lemma}
Let $\mathcal{U}_d \subseteq 
\mathfrak{X}(\K^2)_{\leq d-1}$ be the open and dense set of polynomial vector 
fields of degree $d-1$, 
having exactly $(d-1)^2$ zeros in $\mathcal{P}_{ci} \subset 
{\it Conf} (\K^2, (d-1)^2)$. 
Assume that $\mathcal{P}_{ci}$
has trivial isotropy  group in $\aff(\K^2)$.
In $\mathcal{U}_d$
there exists an analytic $SL(2,\K)$--bundle structure as follows
\begin{equation}\label{fibrado-PSL}
\begin{picture}(200,55)
\put(50,46){$SL(2,\K)$}
\put(94,49){\vector(1,0){35}}
\put(136,46){$\mathcal{U}_d$}
\put(144,41){\vector(0,-1){25}}
\put(151,26){$\pi$}
\put(125,0){$\dfrac{\mathcal{U}_d}{SL(2,\K)} 
\subseteq {\it Conf}(\K^2, (d-1)^2).$}
\end{picture}
\end{equation}

\end{lemma}

\begin{proof}
We want to show that
a polynomial vector field  
$X \in \mathfrak{X}(\K^2)_{\leq d-1}$
has $(d-1)^2$ zeros exactly at $\mathcal{P}_{ci}$ 
as in \eqref{puntos-base-de-pincel} 
if and only if 
it
is of the shape $X_{\tt M}$ in \eqref{pincel-de-campos}.

\noindent
\smallskip
($\Rightarrow$) Let $X=A(x,y)\del{}{x} + B(x,y)\del{}{y}$
be a vector field in $\mathfrak{X}(\K^2)_{\leq d-1}$.
The curve $\mathcal{C}_A \doteq \{ A(x,y)=0\}$ 
has at most degree  $d-1$ and it
would contain $\mathcal{P}_{ci}$.
There exist an open set of values 
$\{ [\mu, \nu ] \} \subset \K \mathbb{P}^1 $ 
such that for each value 
the respective curve  
$\{\mu F + \nu G =0 \}$ in  the pencil \eqref{pincel-curvas}
intersects in a transversal way $\mathcal{C}_A$ at every point of $\mathcal{P}_{ci}$.
By B\'ezout's theorem, 
the degree of $\mathcal{C}$ is exactly $d-1$.  
For any point 
$p \in \mathcal{C}_A \backslash \mathcal{P}_{ci} \subset \K^2$, 
there exists 
a value, say $[-{\tt c} ,  -{\tt d} ]$ in \eqref{pincel-curvas}
such that its respective curve satisfies
$\mathcal{C}_{-{\tt c}  -{\tt d} }  
\cap \mathcal{C}_A \supset \widehat{\mathcal{P}} \cup \{p\}$. 
Hence (again by B\'ezout's theorem), both curves coincide as sets
and $A= -{\tt c}F - {\tt d}  G$ as polynomials.
\end{proof}

Thus, each configuration $\mathcal{P}_{ci}$
has an associated fiber
$
\big\{
X_{\tt M}  \ \vert \ {\tt M} \in SL(2, \K)
\big\} \subset \mathcal{U}_d
$
in \eqref{fibrado-PSL}, which is a family of 
not necessarily Hamiltonian vector fields.
A further goal is the study of the intersection
$$
\big\{ X_{\tt M}  \ \vert \ {\tt M} \in SL(2, \K)
\big\}
\cap Ham(\K^2)_{\leq d}.
$$

\begin{corollary} A jump phenomena.
Let $\mathcal{P} =
\{ (0,0), (1,0), (1/2, \sqrt{3}/2), (x_4, y_4) \}$ 
be a configuration 
leading  to a family of vector fields
$\mathfrak{F} (\mathcal{P})
= \{ X_{\tt m }  \ \vert \ {\tt m} \in SL(2, \K)  \} 
$ 
as in \eqref{pincel-de-campos}.

\begin{enumerate}
\item  
If $(x_4, y_4) \in \K^2 \backslash  {\tt A}$, then 
there exists one projective class in 
$\mathfrak{F}(\mathcal{P}) \cap Ham(\K^2)\leq_2$. 

\item
If $(x_4, y_4) = {\tt R}_1, \,  {\tt R}_2
$ or  ${\tt R}_3$, 
then 
there exists a $\mathbb{KP}^1$--family of Hamiltonian vector fields
$\mathfrak{F}(\mathcal{P})\cap Ham(\K^2)_{\leq_2}$.
\hfill
$\Box$
\end{enumerate} 
\end{corollary}

\begin{example}{\it 
A family $\big\{ X_{\tt M}  \ \vert \ {\tt M} \in SL(2, \K)
\big\}$  in  
\eqref{fibrado-PSL}
with $(d-1)^2 \geq 4$ points as base locus  and
such that
its Hamiltonian vector fields 
$ Ham(\K^2)_{\leq d-1}= [f]$ 
determine one projective class.} 
\begin{upshape}

\noindent 
Consider two algebraic curves such that
$$
\mathcal{P}_{ci}=
\{ 
\underbrace{y- \mu \Pi_{\iota=1}^d (x-x_\iota)=0}_{
F(x,y)=0}
\}
\cap 
\{ 
\underbrace{x- \nu \Pi_{j=1}^d (y- y_j )=0}_{G(x,y)=0} \},
\ \ \ 
d \geq 3
$$

\noindent 
has exactly $(d-1)^2 \geq 4$ points.

It follows  that, the associated 1--form 
$\omega_{\tt m}$ is exact if and only if 
$
{\tt m}=
\left(
\begin{smallmatrix}
\tt{a} & \tt{ 0} \\
\ \tt{0} & \tt{a} 
\end{smallmatrix}
\right) 
$. 
\noindent 
In fact, 
suppose $f(x,y)$ such that 
$\omega_{\tt m} =df$, then  

\centerline{ 
${\tt a}F(x,y)+{\tt b}G(x,y)=f_x$ 
\ and \ 
${\tt c}F(x,y)+{\tt d}G(x,y)=f_y$.} 

\noindent 
As $f_{xy}=f_{yx}$, 
then ${\tt a}- {\tt b} \frac{\partial}{\partial y} \Pi_{j=1}^d (y- y_j )
=-{\tt c}\frac{\partial}{\partial x}\Pi_{\iota=1}^d (x-x_\iota)+{\tt d}$,
so ${\tt a}={\tt d}$ and ${\tt b}={\tt c}=0$.

\noindent 
By assumption $\omega_{\tt m }$ is exact and
defining  $f_{\tt m} (x,y) = \int^{(x,y)} \omega_{\tt m}$,
we conclude that
\begin{equation}
\mathfrak{F}(\mathcal{P}_{ci}) \cap 
Ham(\K^2)_{\leq d-1} =
\mathcal{L}_{d} (\mathcal{P}_{ci})= [f_{\tt m}] 
\ \ \
\hbox{ and }
\ \ \
dim_\K (\mathcal{L}_{d} (\mathcal{P}_{ci} )) = 0 .
\end{equation}
\end{upshape}
\end{example}

\begin{example}{\it 
A fiber $\big\{ X_{\tt M}  \ \vert \ {\tt M} \in SL(2, \K)
\big\}$ as in  \eqref{fibrado-PSL}
with $(d-1)^2\geq 9$ points as a base locus  and
such that
 
\centerline{
$\big\{ X_{\tt M}  \ \vert \ {\tt M} \in SL(2, \K)
\big\} 
\cap Ham(\K^2)_{=d} = \emptyset$.}
}
\begin{upshape}

\noindent 
Consider two hyperelliptic curves such that 
$$
\widehat{\mathcal{P}}=
\{ F(x,y)=y^2- \mu \Pi_{\iota = 1}^d (x-x_\iota )=0 \} 
\cap 
\{ G(x,y)= x^2 - \nu \Pi_{j=1}^d (y-y_j)=0 \}
$$
has  exactly $(d-1)^2 \geq 9$ points.
It follows  that 
$\omega_{\tt m}$ is non exact for all
$
{\tt m}=
\left(
\begin{smallmatrix}
\tt{a} & \tt{ b} \\
\tt{c} & \tt{d} 
\end{smallmatrix}
\right) . 
$
We conclude that
\begin{equation}  
\mathcal{L}_{d} (\widehat{\mathcal{P}})= \emptyset 
\ \ \
\hbox{ and }
\ \ \
dim_\K (\mathcal{L}_{d} ( \widehat{\mathcal{P}})) = -1 .
\end{equation}

\noindent 
In fact,  if we
suppose $f(x,y)$ such that $\omega_m=df$, 
then $2{\tt a}y- {\tt b}\frac{\partial}{\partial y} 
\Pi_{j=1}^d (y- y_j )=-{\tt c}\frac{\partial}{\partial x}
\Pi_{\iota=1}^d (x-x_\iota)+2 {\tt d} x$, 
so ${\tt a}={\tt b}={\tt c}={\tt d}=0$.
\end{upshape}
\end{example}

\begin{corollary}\label{familias-no-Morse}  
There exists a fiber $\mathfrak{F}$ 
as in \eqref{fibrado-PSL} 
having $d^2$ points as base locus and

\centerline{$
\mathfrak{F}(\widehat{\mathcal{P}}) \cap Ham(\K^2)_{=d} = \K \mathbb{P}^1 .  $}

\noindent 
Moreover, $\K \mathbb{P}^1 $ minus a finite set 
determines Morse polynomials.
\end{corollary}

The above result uses the following very particular 
configurations.

\begin{definition}\label{red-d2-puntos}
\begin{upshape}
{\it 
A grid of $(d-1)^2$ points} $\mathcal{G}$
is determined by two sets of $d-1$ parallel lines where 
one set is transverse to the other, {\it i.e.}
up to affine transformation 
$$
\mathcal{G}=
\{ F(x,y)= \Pi_{j=1}^{d-1} (y- y_j )=0 \}
\cap
\{ G(x,y)=\Pi_{\iota=1}^{d-1} (x-x_\iota)=0\}  
$$
with exactly $(d-1)^2 \geq 4 $ points (it is a complete intersection).
\end{upshape}
\end{definition}

\noindent 
{\it Proof of the Corollary.}
The family $X_{\tt M}$ with a grid of $(d-1)^2$ points
is  Hamiltonian if and only if 

\centerline{
$
{\tt M} \in 
\left\{
\left(
\begin{smallmatrix}
0& -{\tt d} &  \\
{\tt a} & 0
\end{smallmatrix}
\right) \right\} \cong \K^2 \subset SL(2, \K)  
$.}

\noindent 
In fact, 
$\omega_{\tt m} =
 ({\tt a} F(x) + {\tt b}G(y)) dx + 
( {\tt c}F(x) +  {\tt d}G(y)) dy =0$ is exact 
if and only if ${\tt b} G(y)_y = {\tt c} F(x)_x$. 
The equality holds only for ${\tt b}= {\tt c}=0$.
   
\noindent 
The respective vector subspace of polynomials 
\begin{equation}\label{polinomios-en-red}
\left\{ f({\tt a}, {\tt d}, x,y) = 
{\tt a} \int^{(x,y)} \Pi_{\iota=1}^d (x-x_\iota) dx + 
{\tt d} \int^{(x,y)} \Pi_{j=1}^d (y-y_j) dy    
\ \Big\vert \
({\tt a}, {\tt d} ) \in \K^2 \backslash \{\overline{0} \} \ \right\} 
\subset  
\K[x,y]^0_{\leq d}
\end{equation}
shows that
\begin{equation}
\mathcal{L}_{d} (\mathcal{P}) \supset \{ [f({\tt a} ,{\tt d}, x,y) ] \} 
\ \ \
\hbox{ and }
\ \ \
dim_\K (\mathcal{L}_{d} (\mathcal{P})) = 1 .
\end{equation}

For $({\tt a}, {\tt d}) \neq ({\tt a}, 0),  \, (0, {\tt d}) $, 
each polynomial $f({\tt a} ,{\tt d}, x,y) \in \K[x,y]^0_{\leq d}$ 
in  \eqref{polinomios-en-red} 
has 
$(d-1)^2$ Morse critical points. In fact, 
at each point $ p\in \mathcal{P}$
a very simple observation with the Taylor series shows that 
$f({\tt a} ,{\tt d}, x,y) = \widetilde{\tt a} x^2 + 
\widetilde{\tt b} y^2 + \mathcal{O}_3(x,y)$,
where $\widetilde{\tt a} \widetilde{\tt b} \neq 0$.

On the other hand, for 
$({\tt a}, {\tt d}) = ({\tt a}, 0), \, (0, {\tt d}) $
the polynomial  $f({\tt a}, {\tt d}, x, y) $ has lines 
of critical points in $\{ P(x,y)=0\}$ or
$\{ Q(x,y)= 0\}$. 
\hfill
$\Box$

\begin{example}\label{familia-de-rotaciones-texto}
Real rotated Hamiltonian vector fields for the grid of $4$ points.   
\begin{upshape} 
Let
$\mathcal{G} = \{ (0,0), (1,0), (0, 1), R=(1,1)\}$ be a grid 
its space of polynomials is

\centerline{$
f({\tt a}, {\tt d}, x,y)=
{\tt a} \Big(\dfrac{x^3}{3}  - \dfrac{x^2}{ 2} \Big) 
+
{\tt d} \Big( \frac{y^3}{3}  - \frac{y^2}{ 2} \Big).
$}

\noindent 
In particular for $\K=\R$, we consider the family

\centerline{
$
R_\theta =\left\{ 
f_\theta (x,y) =
\cos(\theta) \big(\frac{x^3}{3} - \frac{x^2}{2} \big)
+
\sin(\theta) \big(\frac{y^3}{3} - \frac{y^2}{2} \big)  
\ \Big\vert \ 
\theta \in [0, 2\pi]\right\}
$}

\noindent 
of polynomials in \eqref{polinomios-en-red}. 
They originate a family of rotated vector fields, 
see Fig. \ref{atlas-polinomios-cubicos}. 
The algebraic curve
$\{  f_\theta (x,y)+c=0 \}$ is reducible for $\theta=\pi/4$ and $c={1/6}$. 
In this case we get 

\centerline{
$\{ (x+y-1)(2y^2-2xy+2x^2-y-x-1)=0 \} $.
}
\end{upshape}
\end{example}

The following family of vector fields 
is related
to the results in \cite{Ramirez} \S5;
see Fig. \ref{atlas-polinomios-cubicos}, upper row.

\begin{corollary}
\label{espectros}
The one dimensional  holomorphic 
family of Hamiltonian vector fields of the 
polynomials 

\centerline{$
\left\{
f({\tt a}, {\tt d}, x,y)=
{\tt a} \Big(\dfrac{x^3}{3}  - \dfrac{x^2}{ 2} \Big) 
+
{\tt d} \Big( \frac{y^3}{3}  - \frac{y^2}{ 2} \Big)
\ \big\vert \  {\tt ad}=1
\right\} $}

\noindent 
has singularities at 
$\mathcal{G} = \{ (0,0), (1,0), (0, 1), R=(1,1)\}$
and spectra of eigenvalues 

\centerline{$
\big[
[i, \, -i], \, [1,-1], \, [i, \, -i], \, [1,-1] 
\big] $.}
\hfill
$\Box$
\end{corollary}

\begin{corollary}
\label{polinomios-Morse-no-esencialmente-determinados}
For $d \geq 3$, there exist Morse polynomials $f \in \K[x,y]_{=d}^0$
with $(d-1)^2$ singular points 
that are not essentially determined. 
\hfill
$\Box$
\end{corollary}

\section{Closing remarks}
\label{closing-remarks}

Let $\mathfrak{X}(\K^2)_{\leq d-1}$ be
the space of polynomial 
vector fields $\{ X \}$ of at most degree $d-1$ 
on $\K^2$. A general and natural question is
as follows. 
Under what conditions  a polynomial vector field $X$ on $\K^2$ is 
{\it essentially determined}  by its 
configuration of zeros $\mathcal{Z}(X)$ in $\K^2$?

In simple words, a vector field  $X$  
{\it is essentially determined} 
(in  $\mathfrak{X}(\K^2)_{\leq d-1}$) by 
its configuration of zeros
$\mathcal{Z}(X_f)$,

\centerline{
if for any  $Y \in \mathfrak{X}(\K^2)_{\leq d-1} $
satisfying $ \mathcal{Z}(X) \subset \mathcal{Z}(Y)\subset \K^2$, 
then $ X = \lambda Y $. }

\noindent 
Recalling that for affine degree $d$ the number 
of isolated singularities of 
the associated singular holomorphic foliation
$\mathcal{F}(\mathcal{X})$ on the whole $\cp^2$ is $(d-1)^2 + d$,
the hypothesis of multiplicity one 
must be understood for all these points.
Proposition \ref{abierto-cerrado-Zariski} 
confirms that in the Hamiltonian 
case only $\delta(d) \leq (d-1)^2$ points are required.

Recall that X. G\'omez--Mont 
and G. Kempf, \cite{GomezMont-Kempf}, established
in the complex rational case the following deep result, 
that also enlightens the real case.

\smallskip

\noindent 
{\it 
A meromorphic vector field $\mathcal{X}$ on 
$\cp^m$, $m \geq 2$,  of degree $r\geq 2$, with 
critical points having all its zeros of multiplicity one 
is completely determined by its zero set.}

\smallskip

Moreover, J. Artes, J. Llibre, D. Schlomiuk
and N. Vulpe, 
\cite{Artes-Llibre-Vulpe},
\cite{Artes-Llibre-Schlomiuk-Vulpe}
prove the following:

\smallskip 

\noindent 
{\it 
A polynomial vector field $\mathcal{X}$ on $\K^2$ of degree 
two, 
is completely determined by the position of its seven critical points 
(including the points at infinity).}

\smallskip

As far as we know, over $\K=\C$
the more general result is due to A. 
Campillo and J. Olivares,
\cite{Campillo-Olivares}:

\smallskip 

\noindent 
{\it 
A singular holomorphic foliation  
$\mathcal{X}$ on $\cp^2$ of degree $r\geq 2$, w
is completely determined by its singular scheme.}

\smallskip

See C. Alc\'antara {\it et al.} 
\cite{Alcantara-Pantaleon} for recent developments
regarding foliations with multiple points. 
We summarize our results as
follows. 

\begin{corollary}
A polynomial Hamiltonian vector field 
$X_f$ on $\K^2$ of degree two 
is completely determined (in the space of polynomial vector fields of degree $2$, up to a scalar factor 
$\lambda \in \K^*$)
by its zero points, when they are four isolated points 
different from $\{(0,0), \, (1,0), \, (0,1), \, (1,1) \}$,
up to affine transformation. 
\end{corollary}

Our hope is that the explicit 
results in this paper can illustrate
the classification of polynomials $\K[x,y]$
up to algebraic equivalence $Aut(\K^2)$; 
see \cite{Fernandez-de-Bobadilla} and
\cite{Wightwick}
for this order of ideas. 
This potential  application 
is the subject of a future project.


\end{document}